\documentclass{elsarticle}

\usepackage{lineno,hyperref}
\usepackage{amsmath}
\usepackage{amsfonts,amssymb}
\usepackage{esint}
\usepackage{subfigure}
\usepackage{wasysym}
\usepackage{amsthm}
\usepackage{multirow}
\usepackage{algorithm}
\usepackage{algorithmic}

\newtheorem{remark}{Remark}

\modulolinenumbers[5]

\begin{document}

\begin{frontmatter}

\title{A class of high-order weighted compact central schemes for solving hyperbolic conservation laws}


\author[add1]{Hua Shen\corref{mycorrespondingauthor}}
\cortext[mycorrespondingauthor]{Corresponding author}
\ead{huashen@uestc.edu.cn}
\author[add2]{Matteo Parsani}
\ead{matteo.parsani@kaust.edu.sa}
\address[add1]{School of Mathematical Science, University of Electronic Science and Technology of China, Chengdu, Sichuan 611731, China}
\address[add2]{King Abdullah University of Science and Technology (KAUST),
Extreme Computing Research Center (ECRC), Computer, Electrical and Mathematical
Sciences \& Engineering (CEMSE),
Thuwal, 23955-6900, Kingdom of Saudi Arabia}

\begin{abstract}
We propose a class of weighted compact central (WCC) schemes for solving hyperbolic conservation laws.
The linear version can be considered as a high-order extension of the central Lax–Friedrichs (LxF) scheme 
and the central conservation element and solution element (CESE) scheme.
On every cell, the solution is approximated by a $P$th-order polynomial of which all the DOFs are stored and updated separately.
The cell average is updated by a classical finite volume scheme
which is constructed based on space-time staggered meshes such that 
the fluxes are continuous across the interfaces of the adjacent control volumes and,
therefore, the local Riemann problem is bypassed.
The $k$th-order spatial derivatives are updated by a central difference of the $\left(k-1\right)$th-order spatial derivatives at cell vertices.
All the space-time information is calculated by the Cauchy-Kovalewski procedure.
By doing so, the schemes are able to achieve arbitrarily uniform space-time high-order on a super compact stencil with only one explicit time step.
In order to capture discontinuities without spurious oscillations,
a weighted essentially non-oscillatory (WENO) type limiter is tailor-made for the schemes.
The limiter preserves the compactness and high-order accuracy of the schemes.
The accuracy, robustness and efficiency of the schemes are verified by several numerical examples of scalar conservation laws and the compressible Euler equations.
\end{abstract}

\begin{keyword}
compact scheme \sep  central scheme \sep CESE scheme \sep high-order scheme \sep finite volume scheme \sep hyperbolic conservation laws
\end{keyword}

\end{frontmatter}


\section{Introduction}\label{sec:intro}
In continuum physics, many phenomena are governed by hyperbolic conservation laws.
Most of them are nonlinear, thereby obtaining the analytical solutions is difficult, if not impossible.
Therefore, the study of numerical methods for hyperbolic conservation laws attracts substantial contemporary interest.

When solving nonlinear hyperbolic conservation laws, one of the challenges is that discontinuous solutions may appear
even when the initial conditions are sufficiently smooth.
Therefore, shock-capturing (or discontinuity-capturing) is an essential feature of an effective numerical scheme for solving general hyperbolic conservation laws.
The first-order upwind scheme and the Lax–Friedrichs (LxF) central scheme \cite{Lax1954} are forerunners of such schemes.
The prototype of upwind schemes can be traced back to the CIR scheme proposed by Courant, Isaacson, and Rees \cite{Courant1952solution}.
By comparison, upwind schemes need to determine the propagating direction of a wave in advance,
but the LxF scheme does not require the structure of the solution for a specific problem.
In this sense, central schemes are more general schemes for solving hyperbolic conservation laws than upwind schemes.
Godunov \cite{Godunov1959} extended the first-order upwind scheme to solve the Euler equations by incorporating an exact nonlinear Riemann solver into the scheme.
For a general hyperbolic system, it is probably difficult to, or even not possible to get the exact Riemann solution,
thereby limiting the application of the scheme.
Afterwards, researchers realized that the Riemann solver is not necessarily exact,
and a lot of approximate Riemann solvers were developed to simplify Godunov scheme, such as Roe Riemann solver \cite{Roe1981JCP},
HLL Riemann solver \cite{Harten1983HLL}, HLLC Riemann solver \cite{Toro1994HLLC}, and so forth.

The first-order upwind scheme and the LxF scheme have large numerical dissipation.
As a result, they require very fine meshes to capture small structures in smooth regions and to obtain a sharp transition for discontinuities.
Therefore, second- and higher-order methods are desired for real applications, but it is not an easy task.
Godunov \cite{Godunov1959} proved that a linear monotonicity-preserving scheme for the linear advection equation is at most first-order.
It means that, even for the linear advection equation, spurious oscillations near discontinuities are inevitable if we use a high-order linear scheme to solve it.
In order to get an oscillation-free high-order scheme for solving hyperbolic conservation laws, nonlinear operators are necessary.
In general, the nonlinear operators can be categorized into three major approaches.
The first approach is the artificial viscosity technique proposed by von Neumann and Richtmyer \cite{vonNeumann1950ArtificialViscosity}
which is still widely used today.
This technique uses a uniform scheme which is accurate in smooth regions but not oscillation-free near discontinuities.
The key idea is using a solution-dependent nonlinear function to tune the numerical viscosity
that is negligible in smooth regions to preserve accuracy but is large enough near discontinuities to kill spurious oscillations.
The second approach is the flux-corrected transport (FCT) technique which was proposed by Boris and his coworkers
\cite{Boris1973FCTI, Book1975FCTII, Boris1976FCTIII, Zalesak1979FCT2D}.
This approach uses a nonlinear flux limiter to hybridize a first-order flux and a second-order flux.
In smooth regions, the limiter will choose the second-order flux to achieve high accuracy;
while near discontinuities, the limiter will choose the first-order one to tame spurious oscillations.
The third approach is the monotonic upstream-centered scheme for conservation laws (MUSCL) proposed by van Leer
\cite{vanLeer1974MUSCLII, vanLeer1977MUSCLIII, vanLeer1977MUSCLIV, vanLeer1979MUSCLV}.
This approach uses a nonlinear limiter to reconstruct the slope of the solution approximated by a linear polynomial.
The philosophy of MUSCL and FCT is similar.
When it degrades to the linear advection case, MUSCL is actually identical with FCT.
It is interesting to note that the Russian scientist Kolgan independently proposed a similar approach as MUSCL \cite{Kolgan1972, vanLeer2011historical}.
In principle, the above mentioned approaches can apply to both upwind schemes and central schemes,
but most of the early researches are second-order extensions of the first-order upwind scheme including the original MUSCL and FCT.
In order to bypass the local Riemann problem at cell interfaces,
there are some high-resolution Riemann-solver-free central schemes extended from the first-order LxF scheme.
The Nessyahu-Tadmor (N-T) scheme \cite{Nessyahu1990NT} and the central space-time conservation element and solution element (CESE) scheme
\cite{Chang1991NASATM, Chang1995JCP} are precursors of such schemes.
The N-T scheme and the CESE scheme are fundamentally similar.
Both of them are second-order schemes constructed on space-time staggered mesh and adopt a slope limiter to achieve the non-oscillatory feature.
The major difference lies in the treatment of the slope.
The N-T scheme first updates the cell averages and then reconstructs the slope using cell averages like the MUSCL does.
However, the CESE scheme considers the slope as an additional degree of freedom (DOF) and updates it by a weighted central difference scheme.

The above mentioned second-order schemes paved the way for constructing high-order non-oscillatory schemes.
Colella and Woodward \cite{Colella1984PPM} further extended MUSCL to the third-order piecewise parabolic method (PPM).
Harten \cite{Harten1983TVD} introduced the concept of total variation (TV)
and proved that MUSCL is total variation diminishing (TVD) in the scalar case \cite{Harten1987ENO}.
By realising that a TVD scheme, regardless its particular form, is only first-order at local extrema.
Harten and his coworkers \cite{Harten1987ENO} slightly relaxed the TVD constraint and constructed uniformly high-order accurate essentially non-oscillatory (ENO) schemes. ENO schemes achieve the non-oscillatory feature by the strategy of choosing the smoothest stencil from several local candidates.
In the light of ENO schemes, Liu \emph{et al}. \cite{Liu1994WENO} proposed weighted ENO (WENO) schemes by assigning proper weights to all the candidate stencils.
Jiang and Shu \cite{Jiang_Shu1996WENO} and Borges \emph{et al}. \cite{Borges2008WENO_Z} further improved the accuracy of WENO schemes by carefully design the weights.
Bianco \emph{et al}. \cite{Bianco1999CENO} and Levy \emph{et al}. \cite{Levy1999CWENO, Levy2000CWENO_ANM} extended the LxF central scheme to high-order
by utilizing ENO and WENO reconstruction respectively.
In principle, ENO and WENO schemes can be extended to arbitrary order, but the stencil contains more and more cells as the order increases.
This feature is not friendly for implementing ENO and WENO schemes on unstructured meshes, although we can do it,
see for example the work of Abgrall \cite{Abgrall1994ENO_Tri} and the work of Hu and Shu \cite{Hu1999WENO_Tri}.
In order to reduce the size of the stencil, Qiu and Shu \cite{Qiu2004HermiteWENO, Qiu2005HermiteWENO2D} proposed a fifth-order Hermite WENO scheme
which utilizes both the function and its first-order derivative in the reconstruction and has a more compact stencil than the original fifth-order WENO scheme.
Levy \emph{et al}. \cite{Levy2000CompactCWENO} proposed a third-order compact central WENO scheme by combining a second-order polynomial with several first-order polynomials. Nevertheless, both the Hermite WENO scheme and the compact central WENO scheme have to enlarge the stencil for higher-order extensions.
There are plenty of other variations of WENO schemes which we are not going to review here.
One major strategy to build compact high-order schemes is increasing the internal DOF,
such as the discontinuous Galerkin (DG) scheme \cite{Reed1973DG, Cockburn1989RKDGII, Cockburn1989RKDGIII, Cockburn1990RKDGIV, Cockburn1998RKDGV},
the compact least-squares finite-volume (CLSFV) scheme \cite{Wang2016CompactFVI, Wang2016CompactFVII, Wang2017CompactFVIII},
the high-order CESE scheme \cite{Liu2004HighOrderCESE, Shen2015HybridCESE, Bilyeu2014HighOrderCESE},
the spectral volume (SV) scheme \cite{Wang2002SV, Wang2002SVII, Wang2004SVIII},
the spectral difference (SD) scheme \cite{Wang2006SDI, Wang2007SDII}.
The DG, CLSFV and CESE schemes increase the internal DOF by increasing the degree of the polynomial,
while SV and SD schemes increase the internal DOF by increasing the number of sub-cells and sub-points respectively.
Based on the finite volume scheme and the DG scheme, Dumbser \emph{et al}. \cite{Dumbser2008PNPM, Dumbser2009PNPM, Dumbser2010PNPM}
proposed a general framework to build high-order $P_NP_M$ schemes, where $N$ is the internal DOF and $M$ is the reconstructed DOF.
The $P_NP_M$ scheme can flexibly choose the internal DOF, thereby adjusting the compactness of the stencil.

As for the time discretization, there are two major strategies.
The first strategy is the method of lines which first discretizes the space and leaves a system of ordinary differential equations (ODEs)
to which any numerical method for initial value problems of ODEs can be applied.
One of the most popular high-order methods is the TVD Runge–Kutta (RK) time discretization
which is simply a convex combination of several first-order Euler time discretization \cite{Shu1988EfficientENO, Jiang_Shu1996WENO, Cockburn1989RKDGII, Cockburn1989RKDGIII, Cockburn1990RKDGIV, Cockburn1998RKDGV}.
Another strategy is the Lax-Wendroff-type time discretization \cite{LaxWendroff1960} which is also called Cauchy-Kovalewski procedure.
This approach repeatedly uses the original partial differential equations (PDEs) to calculate all the space-time information
so that the scheme can achieve high-order accuracy within only one time-step.
Because of the one-time-step feature, the Lax-Wendroff-type time discretization has been used in many schemes,
such as the original ENO scheme \cite{Harten1987ENO}, the WENO finite difference scheme with Lax-Wendroff-type time discretization  \cite{Qiu2003LaxW_WENOFD},
the DG scheme with Lax-Wendroff-type time discretization \cite{Qiu2005LaxWDG}, the $P_NP_M$ scheme \cite{Dumbser2008PNPM, Dumbser2009PNPM, Dumbser2010PNPM},
the ADER (arbitrary high-order schemes utilizing high-order derivatives) scheme \cite{Toro2002ADER, Toro2005ADER, Toro2006ADER, Balsara2013ADER},
the CESE scheme \cite{Chang1991NASATM, Chang1995JCP, Bilyeu2014HighOrderCESE, Liu2004HighOrderCESE, Shen2015HybridCESE}, and so forth.

Second-order schemes significantly improve the resolution and narrow the thickness of discontinuities as comparing with first-order schemes.
In academia and engineering, most people reach a consensus that second-order schemes are more efficient than first-order schemes for most problems.
However, it seems that there is some controversy associated with the efficiency of second-order schemes and high-order (higher than second-order) schemes.
In order to achieve a certain level of accuracy, we can either use a second-order scheme on a fine mesh or use a high-order scheme on a relatively coarse mesh.
Although, on the same mesh, a high-order scheme is indeed more expensive than a second-order scheme based on the same framework,
a high-order scheme converges much faster than a second-order scheme
because the error of a $P$th-order scheme is proportional to $\Delta x^P$ where $\Delta x$ is the mesh size.
As a result, a high-order scheme is more efficient than a second-order scheme to achieve a high-level of accuracy \cite{Wang2007PAS, Wang2013IJNMF}.
In many applications, we have to simultaneously solve discontinuities and highly sophisticated structures, such as detonation waves \cite{Henrick2006Det}, interface instabilities induced by a shock impact \cite{Latini2007RMI}, shock-turbulence interaction \cite{Johnsen2010JCP}, and so forth.
This is the major motivation to develop high-order shock-capturing schemes.

In this work, we first introduce a class of arbitrarily high-order linear schemes for solving hyperbolic conservation laws.
The schemes can be considered as an extension of the central conservation element and solution element (CESE) scheme \cite{Chang1995JCP, Liu2004HighOrderCESE, Shen2015HybridCESE}.
On every cell, the solution is approximated by a $P$th-order polynomial of which all the DOFs are stored and updated separately.
The cell average is updated by a classical finite volume scheme
which is constructed based on space-time staggered meshes such that 
the fluxes are continuous across the interfaces of the adjacent control volumes and, 
therefore, a local Riemann solver is not required.
The $k$th-order spatial derivatives are updated by a central difference of the $\left(k-1\right)$th-order spatial derivatives at cell vertices.
All the space-time information is calculated by the Cauchy-Kovalewski procedure.
In this way, the schemes are able to achieve uniformly space-time high-order on a super compact stencil with only one explicit time-step.
In addition, a compact WENO-type limiter is designed to capture discontinuities.
Because of the above features of the schemes, we call them weighted compact central (WCC) schemes for solving hyperbolic conservation laws.

\section{Arbitrarily high-order WCC schemes in one dimension}\label{SEC:1DWCCS}
\subsection{Linear compact central schemes in one dimension}\label{SubSEC:1DLinearWCCS}
For the convenience of expression, we consider the one-dimensional scalar conservation law, equipped with certain initial conditions and boundary conditions,
which can be expressed in the following form
\begin{subequations}\label{Eq:1DHCL}
  \begin{equation}\label{SubEq:1DHCL_Eq}
    \frac{\partial u}{\partial t}+\frac{\partial f(u)}{\partial x}=0,\quad x\in[x_L,x_R],\quad t\in[0,\infty),
  \end{equation}
  \begin{equation}\label{SubEq:1DHCL_IC}
    u(x,0)=u_0(x), \quad x\in[x_L,x_R],
  \end{equation}
  \begin{equation}\label{SubEq:1DHCL_BC}
    u(x_L,t)=u_L(t),\quad u(x_R,t)=u_R(t), \quad t\in[0,\infty).
  \end{equation}
\end{subequations}

\begin{figure}
  \centering
  \includegraphics[width=7 cm]{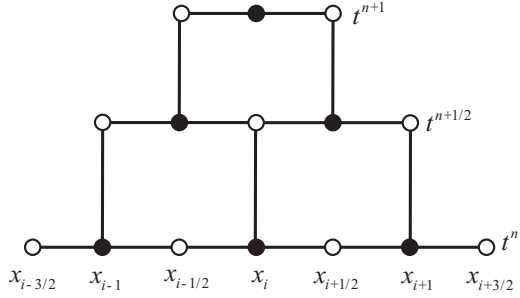}
  \caption{Schematic diagram of the space-time staggered mesh for 1D WCC schemes.
           Cell vertices are denoted by $\Circle$ and cell centers are denoted by $\CIRCLE$.}
  \label{FIG:1DMesh}
\end{figure}

The spatial domain $[x_L,x_R]$ is discretized by two mesh systems.
One mesh is defined by $x_i=x_L+(i-1)\Delta x$ ($i=1 \text{ to } N+1$), where $\Delta x=(x_R-x_L)/N$.
The other mesh is defined by $x_{i+1/2}=(x_i+x_{i+1})/2$ ($i=0 \text{ to } N+1$),
where $x_{1/2}$ and $x_{N+3/2}$ are ghost mesh points.
The solution of Eq. (\ref{Eq:1DHCL}) is updated by a staggered strategy as shown in Fig. \ref{FIG:1DMesh}:
\begin{itemize}
  \item At $t=t^n$, the solution is stored on the original cell $\mathcal{C}_i$ that is defined by $[x_{i-1/2},x_{i+1/2}]$.
  \item At $t=t^{n+1/2}$, the solution is shifted to the staggered cell $\mathcal{C}_{i+1/2}$ that is defined by $[x_{i},x_{i+1}]$.
  \item At $t=t^{n+1}$, the solution reverts to $\mathcal{C}_i$, preparing for the next loop.
\end{itemize}

For the seek of concise expression, we denote a scaled arbitrary-order partial derivative of a sufficiently continuous function $\psi(x_1, x_2, ...)$ as
\begin{equation}\label{Eq:1D_Derivative_DEF}
  \psi_{k_1x_1,k_2x_2,...}=\frac{\partial^{(k_1+k_2+\cdots)} \psi}{\partial x_1^{k_1}\partial x_2^{k_2}\cdots}\Delta x_1^{k_1}\Delta x_2^{k_2}\cdots.
\end{equation}
On cell $\mathcal{C}_i$ at $t^{n}$, the solution is approximated by a $P$th-order polynomial as
\begin{equation}\label{Eq:1DPolynomial_Space}
  u(x,t^{n})=\sum_{k=0}^{P} \frac{1}{k!}(u_{kx})_i^{n}\left(\frac{x-x_i}{\Delta x}\right)^k,\quad x\in \mathcal{C}_i,
\end{equation}
where $(u_{kx})_i^{n}$ is the scaled $k$th-order spatial derivatives of the approximated solution at $(x_i,t^n)$.

The time marching scheme from $t=t^n$ to $t=t^{n+1/2}$ is almost the same as that from $t=t^{n+1/2}$ to $t=t^{n+1}$,
except that the mesh is shifted by $\Delta x/2$.
Without loss of generality, we consider the time marching scheme from $t=t^n$ to $t=t^{n+1/2}$.
First, we perform the integral of Eq. (\ref{Eq:1DHCL}) on the control volume $[x_{i},x_{i+1}]\times[t^n,t^{n+1/2}]$,
obtaining the cell average by a classical finite volume scheme as
\begin{subequations}\label{Eq:1DIntegral}
  \begin{equation}\label{SubEq:1DIntegral}
    \bar{u}_{i+1/2}^{n+1/2}=\frac{1}{2}[(\bar{u}_R)_i^n+(\bar{u}_L)_{i+1}^n]+\frac{\Delta t}{\Delta x}(\bar{f}_{i}^n-\bar{f}_{i+1}^n),
  \end{equation}
  \text{where}
  \begin{equation}\label{SubEq:1DCellAVG}
    \bar{u}_{i+1/2}^{n+1/2}=\frac{1}{\Delta x}\int_{x_{i}}^{x_{i+1}}u(x,t^{n+1/2})dx
    =\sum_{k=0}^{P} \frac{1+(-1)^k}{(k+1)!2^{(k+1)}}(u_{kx})_{i+1/2}^{n+1/2},
  \end{equation}
  \begin{equation}\label{SubEq:1DCellAVGR}
    (\bar{u}_R)_i^{n}=\frac{2}{\Delta x}\int_{x_{i}}^{x_{i+1/2}}u(x,t^{n})dx
    =\sum_{k=0}^{P} \frac{1}{(k+1)!2^k}(u_{kx})_{i}^{n},
  \end{equation}
  \begin{equation}\label{SubEq:1DCellAVGL}
    (\bar{u}_L)_{i+1}^{n}=\frac{2}{\Delta x}\int_{x_{i+1/2}}^{x_{i+1}}u(x,t^{n})dx
    =\sum_{k=0}^{P} \frac{(-1)^k}{(k+1)!2^k}(u_{kx})_{i+1}^{n},
  \end{equation}
  \begin{equation}\label{SubEq:1DFlux}
    \Delta t=t^{n+1/2}-t^n,\quad \bar{f}_{i}^n=\frac{1}{\Delta t}\int_{t^n}^{t^{n+1/2}}f(x_{i},t)dt.
  \end{equation}
\end{subequations}
In order to completely determine $\bar{u}_{i+1/2}^{n+1/2}$, we still need to calculate $\bar{f}_{i}^n$.
As seen from Fig. \ref{FIG:1DMesh}, when $\Delta t$ satisfies the Courant–Friedrichs–Lewy (CFL) condition,
waves emitted from the cell interfaces do not interact with each other.
As a result, any physical quantity $\psi$ across the left and right boundaries of control volume $[x_{i},x_{i+1}]\times[t^n,t^{n+1/2}]$
is still continuous thanks to the staggered mesh.
Therefore, we may approximate $\psi(x_{i},t)$ by a $P$th-order polynomial as
\begin{equation}\label{Eq:1DPolynomial_Time}
  \psi(x_i,t)=\sum_{k=0}^{P} \frac{1}{k!}(\psi_{kt})_i^n\left(\frac{t-t^n}{\Delta t}\right)^k, t\in [t^n,t^{n+1/2}].
\end{equation}
Substituting Eq. (\ref{Eq:1DPolynomial_Time}) into Eq. (\ref{SubEq:1DFlux}), we obtain
\begin{equation}\label{Eq:Cal1DFlux}
    \bar{f}_{i}^{n}=\sum_{k=0}^{P} \frac{1}{(k+1)!}(f_{kt})_i^{n}.
\end{equation}

Next, we use a central difference scheme to update the spatial derivatives $(u_{(k+1)x})_{i+1/2}^{n+1/2}$ $(0\le k<P)$.
The Taylor series on $\mathcal{C}_{i+1/2}$ imply that
\begin{subequations}\label{Eq:Cal_1D_Derivatives}
\begin{equation}\label{SubEq:Cal_1D_Derivatives_1}
  (u_{kx})_i^{n+1/2}=\sum_{m=0}^{P-k} \frac{1}{m!}(u_{(k+m)x})_{i+1/2}^{n+1/2}\left(-\frac{1}{2}\right)^{m},
\end{equation}
\begin{equation}\label{SubEq:Cal_1D_Derivatives_2}
  (u_{kx})_{i+1}^{n+1/2}=\sum_{m=0}^{P-k} \frac{1}{m!}(u_{(k+m)x})_{i+1/2}^{n+1/2}\left(\frac{1}{2}\right)^{m}.
\end{equation}
\end{subequations}
Subtracting Eq. (\ref{SubEq:Cal_1D_Derivatives_1}) from Eq. (\ref{SubEq:Cal_1D_Derivatives_2}), we obtain
\begin{subequations}\label{Eq:Cal_1D_Derivatives_2}
  \begin{equation}\label{SubEq:Cal_1D_Derivatives_3}
   (u_{(k+1)x})_{i+1/2}^{n+1/2}=(u_{kx})_{i+1}^{n+1/2}-(u_{kx})_{i}^{n+1/2}
      -\sum_{m=3}^{P-k} \frac{1-(-1)^{m}}{m!2^{m}}(u_{(k+m)x})_{i+1/2}^{n+1/2},
  \end{equation}
  \text{where}
  \begin{equation}\label{SubEq:Cal_1D_Derivatives_4}
   (u_{kx})_{i}^{n+1/2}=\sum_{m=0}^{P-k} \frac{1}{m!}(u_{kx,mt})_{i}^{n}.
  \end{equation}
\end{subequations}
We note that, when $P-k<3$ the last term of Eq. (\ref{SubEq:Cal_1D_Derivatives_3}) vanishes,
otherwise it contains higher than $(k+1)$th-order derivatives at $t^{n+1/2}$.
In order to get an explicit expression, we have to update the derivatives from high-order to low-order.

As we can see, Eqs. (\ref{Eq:Cal1DFlux}) and (\ref{SubEq:Cal_1D_Derivatives_4}) involve the time derivatives of $f$ and $u$.
They can be calculated by the Cauchy-Kovalewski procedure.
The derivatives of $f$ are calculated by the chain rule in calculus,
and the time derivatives of $u$ are calculated by repeatedly using the original PDE as
\begin{equation}\label{Eq:1DCauchy_Kovalewski}
  \frac{\partial^{a+b} u}{\partial x^a \partial t^b}=-\frac{\partial^{a+b} f}{\partial x^{a+1}\partial t^{b-1}}.
\end{equation}
Then, $\bar{u}_{i+1/2}^{n+1/2}$ can be explicitly updated by Eq. (\ref{SubEq:1DIntegral})
with the aid of Eqs. (\ref{SubEq:1DCellAVGR}), (\ref{SubEq:1DCellAVGL}), and (\ref{Eq:Cal1DFlux}),
and the spatial derivatives can be explicitly updated by Eq. ({\ref{Eq:Cal_1D_Derivatives_2}).
Finally, we rearrange Eq. (\ref{SubEq:1DCellAVG}) and obtain
\begin{equation}\label{Eq:1DCalU}
    u_{i+1/2}^{n+1/2}=\bar{u}_{i+1/2}^{n+1/2}-\sum_{k=2,k\le P}^{P} \frac{1+(-1)^k}{(k+1)!2^{(k+1)}}(u_{kx})_{i+1/2}^{n+1/2}.
\end{equation}
At this point, all the DOFs in $\mathcal{C}_{i+1/2}$ at $t^{n+1/2}$ have been updated.

\begin{remark}
   The proposed compact central schemes with arbitrarily high-order can be formally considered as
   an extension of the central CESE schemes \cite{Chang1995JCP, Liu2004HighOrderCESE, Shen2015HybridCESE}.
   However, the conservation element (CE) and the solution element (SE) defined in the CESE schemes
   are not building blocks of the schemes, so we do not use the terminology `CESE'.
\end{remark}

\begin{remark}
   The compact central schemes can be considered as special finite-volume-type schemes.
   They use a standard finite-volume scheme, Eq. (\ref{SubEq:1DIntegral}), to update cell averages, so they are naturally conservative.
   The main differences between the present schemes and standard finite-volume schemes are the calculation of spatial derivatives
   and the time-marching method.
   Standard finite-volume schemes reconstruct spatial derivatives using cell averages and require a large stencil to achieve high-order in space.
   In contrast, the present schemes use a central difference within the local control volume to calculate spatial derivatives, thereby having a super compact stencil, as shown by Eq. (\ref{Eq:Cal_1D_Derivatives_2}).
   Regarding the time stepping, standard finite-volume schemes usually use the method of lines which first discretizes the space,
   and then use a multi-step RK method to achieve high-order in time.
   However, the present schemes use the Cauchy-Kovalewski procedure to calculate all the space-time information and achieve
   the same order in time and space by using only one time step, which is similar to the ADER schemes \cite{Toro2002ADER, Toro2005ADER, Toro2006ADER, Balsara2013ADER}.
\end{remark}

\begin{remark}
   The present schemes are purely central schemes based on two facts: (i) the update of cell averages is Riemann-solver-free;
   and (ii) spatial derivatives are updated by central difference schemes.
\end{remark}
\subsection{Weighted compact central schemes in one dimension}\label{SubSEC:1DWCC}
Regardless of particular forms, high-order linear schemes always produce spurious oscillations near discontinuities.
In order to tame the spurious oscillations, we design a non-linear WENO-type limiter for the high-order compact central schemes.

We denote the original polynomial Eq. (\ref{Eq:1DPolynomial_Space}) derived
by the $(P+1)$th-order linear compact central scheme as
\begin{subequations}\label{Eq:1D_U012}
  \begin{equation}\label{SubEq:1DU0}
  u_0^n(x)=\sum_{k=0}^{P} \frac{1}{k!}(u_{kx}^0)_{i}^n\left(\frac{x-x_i}{\Delta x}\right)^k,\quad x\in \mathcal{C}_i.
\end{equation}
\text{Additionally, we construct another two first-order polynomials}
\begin{equation}\label{SubEq:1DU12}
  u_m^n(x)=(u^m)_i^n+(u_x^m)_{i}^n\left(\frac{x-x_i}{\Delta x}\right),\quad x\in \mathcal{C}_i,\quad m=1,2,
\end{equation}
\text{where}
\begin{equation*}
  (u^1)_i^n=(u^2)_i^n=\overline{u_0^n(x)}, (u_x^1)_{i}^n=2\left((u^1)_i^n-u_{i-1/2}^n\right), (u_x^2)_{i}^n=2\left(u_{i+1/2}^n-(u^2)_i^n\right),
\end{equation*}
\end{subequations}
$\overline{u_0^n(x)}$ is the cell average of $u_0^n(x)$.
We note that, the construction of $u_1^n(x)$ and $u_2^n(x)$ is rather straightforward,
since $\overline{u_0^n(x)}$ and $u_{i\pm1/2}^n$ are respectively derived by Eq. (\ref{SubEq:1DIntegral}) and Eq. (\ref{SubEq:Cal_1D_Derivatives_4}) already.

The final limited polynomial is defined as
\begin{equation}\label{Eq:1DLimitedPolynomial}
  u_*^n(x)=\sum_{m=0}^{2}w_mu_m^n(x) \text{ with }w_m=\frac{\hat{w}_m}{\hat{w}_0+\hat{w}_1+\hat{w}_2}.
\end{equation}
Apparently, $\sum_{m=0}^{2}w_m=1$, so it is easy to verify that $\bar{u}_*^n=\bar{u}_0^n$
which means the limited polynomial preserves conservativeness.
Depending on two different cases, the weights are calculated in the following two ways:
\begin{enumerate}
  \item In the case of the second-order scheme ($P=1$), $u_m^n(x)$ $(m=0,1,2)$ have a uniform degree.
  Following Jiang and Shu \cite{Jiang_Shu1996WENO}, the weights are calculated as
  \begin{equation}\label{Eq:1D2ndOrderWeights}
    \hat{w}_m=\left(\frac{1}{\beta_m+\epsilon}\right)^\alpha,
  \end{equation}
  where $\epsilon=10^{-40}$ is used to avoid dividing by 0,
  and $\beta_m$ is the smoothness indicator of $u_m^n(x)$.
  The smoothness indicator of a $P$th-order polynomial defined on $\mathcal{C}_i$ is calculated as \cite{Jiang_Shu1996WENO}
  \begin{equation}\label{Eq:1DSmoothIndicator}
    \beta_m=\sum_{k=1}^{P}\Delta x^{2k-1}\int_{x_i-\Delta x/2}^{x_i+\Delta x/2} \left(\frac{d^ku_m^n(x)}{dx^k}\right)^2dx.
  \end{equation}
  In this case, the limiter is a TVD-like limiter which degrades a second-order scheme to first-order at local extrema \cite{Harten1987ENO}.

  \item In the case of the higher than second-order schemes ($P\ge2$), $u_m^n(x)$ $(m=0,1,2)$ have unequal degrees.
  The weights are modified from Borges \emph{et al}. \cite{Borges2008WENO_Z}
  \begin{subequations}\label{Eq:1DHighOrderWeights}
    \begin{equation}
    \hat{w}_0=1+\left(\frac{\sigma\tau}{\beta_0^2+\epsilon}\right)^\alpha,\quad \hat{w}_m=\left(\frac{\sigma\tau}{\beta_m^2+\epsilon}\right)^\alpha (m=1,2),
   \end{equation}
   \text{where}
    \begin{equation}
    \sigma=\left[u_0^n\left(x_{i+1/2}\right)-u_{i+1/2}^n\right]^2+\left[u_0^n\left(x_{i-1/2}\right)-u_{i-1/2}^n\right]^2,
   \end{equation}
   \begin{equation}
    \tau=\Delta x^{2P-1}\int_{x_i-\Delta x/2}^{x_i+\Delta x/2} \left(\frac{d^Pu_0^n(x)}{dx^P}\right)^2dx=\left[(u_{Px}^0)_{i}^n\right]^2.
   \end{equation}
  \end{subequations}
   We note that $u_0^n\left(x_{i\pm1/2}\right)$ is derived by a space approximation using Eq. (\ref{SubEq:1DU0}),
   while $u_{i\pm1/2}^n$ is derived by a time approximation using Eq. (\ref{SubEq:Cal_1D_Derivatives_4}).
   Therefore, $\sigma$ is the measurement of the residuals at $x_{i\pm 1/2}$,
   and $\tau$ is the smoothness indicator for the highest-order term of $u_0^n(x)$.
   For a $(P+1)$th-order scheme in smooth region, we have $\sigma\propto O\left(\Delta x^{2(P+1)}\right)$, $\tau\propto O\left(\Delta x^{2P}\right)$
   and
   \begin{equation}
    \beta_m\propto\begin{cases}
                    O\left(\Delta x^2\right), & \mbox{if } (u_x)_i^n\ne 0,\\
                    O\left(\Delta x^4\right), & \mbox{if } (u_x)_i^n=0 \mbox{ and } (u_{2x})_i^n\ne 0,
                  \end{cases}\quad (m=0,1,2),
   \end{equation}
   which imply that
   \begin{equation}\label{Eq:1DWmError}
    w_m\propto\begin{cases}
                    O\left(\Delta x^{\alpha (4P-2)}\right), & \mbox{if } (u_x)_i^n\ne 0\\
                    O\left(\Delta x^{\alpha (4P-6)}\right), & \mbox{if } (u_x)_i^n=0 \mbox{ and } (u_{2x})_i^n\ne 0,
                  \end{cases}\quad (m=1,2).
   \end{equation}
   Since $\sum_{m=0}^{2}w_m=1$, the difference between the limited polynomial and the original solution can be written as
   \begin{equation}\label{Eq:1DPolynomialError}
    \varepsilon=u_*^n(x)-u_0^n(x)=\sum_{m=1}^{2}w_m(u_m^n(x)-u_0^n(x)).
   \end{equation}
   Since $u_0^n(x)$ is a $(P+1)$th-order ($P\ge 2$) solution and $u_1^n(x)$ and $u_2^n(x)$ are 2nd-order solutions,
   we have $\varepsilon\propto O\left(\Delta x^2\right)$.
   Combining Eqs. (\ref{Eq:1DWmError}) and (\ref{Eq:1DPolynomialError}), we have the following estimation
   \begin{equation}
    \varepsilon\propto\begin{cases}
    O\left(\Delta x^{\alpha (4P-2)+2}\right), & \mbox{if } (u_x)_i^n\ne 0\\
    O\left(\Delta x^{\alpha (4P-6)+2}\right), & \mbox{if } (u_x)_i^n=0 \mbox{ and } (u_{2x})_i^n\ne 0.
    \end{cases}
   \end{equation}
   In order not to reduce the order of the original solution, it is sufficient to guarantee
   \begin{equation}
    \alpha (4P-6)+2\ge P+1\Rightarrow \alpha \ge\frac{P-1}{4P-6}=\frac{1}{4}+\frac{1}{8P-12}, \quad P\ge2.
   \end{equation}
   When $P=2$, the right hand side of the inequality achieves the maximum value $1/2$,
   so it is sufficient to set $\alpha\ge 1/2$.

   At discontinuities, the high-order solution $u_0^n(x)$ is not smooth,
   and we have $\sigma\propto O\left(1\right)$, $\tau\propto O\left(1\right)$, and $\beta_0\propto O\left(1\right)$.
   However, one of the first-order polynomial $u_m^n(x)$ ($m$=1 or 2) is still smooth,
   and we have $\beta_m\propto O\left(\Delta x^2\right)$ ($m$=1 or 2).
   Therefore, the weight of $u_0^n(x)$ becomes very small, and $u_*^n(x)$ approaches to the smooth first-order polynomial.
   In this case, the weighted procedure behaves like a TVD limiter
   of which the strength increases with the increasing of $\alpha$.
   In practice, we set $\alpha=2$.
\end{enumerate}

\begin{remark}
The above high-order limiter is a variation of the WENO schemes with adaptive order that combine polynomials with different degrees as
\cite{Levy2000CompactCWENO, Zhu2016NewWENOFD, Balsara2016WENO_AO, Zhu2017NewWENOFV_Tri, Zhu2018NewMultiResolutionWENO, Zhu2019NewMultiResolutionWENO_Tri, Zhu2020NewMultiResolutionWENO_Tet}
\begin{equation}\label{Eq:1DLimitedPolynomialOld}
  u_*^n(x)=w_0[\frac{1}{\gamma_0}u_0^n(x)-\frac{\gamma_1}{\gamma_0}u_1^n(x)-\frac{\gamma_2}{\gamma_0}u_2^n(x)]+w_1u_1^n(x)+w_2u_2^n(x),
\end{equation}
where $\gamma_m$ (m=0, 1, 2) are linear weights.
Obviously, the proposed limiter, Eq. (\ref{Eq:1DLimitedPolynomial}),  is more concise than Eq. (\ref{Eq:1DLimitedPolynomialOld}).
Moreover, Eq. (\ref{Eq:1DLimitedPolynomial}) is a convex combination of three solutions,
but the first term on the right hand side of Eq. (\ref{Eq:1DLimitedPolynomialOld}) contains negative coefficients of $u_1^n(x)$ and $u_2^n(x)$
which implies that Eq. (\ref{Eq:1DLimitedPolynomialOld}) is not a strictly convex combination of the three solutions.
\end{remark}

\subsection{Linear stability analysis}\label{SubSEC:1DLinearStability}
In order to provide a guidance for determining the time-step size,
we perform von Neumann stability analysis to assess the linear stability of the proposed WCC schemes.
As shown by Eqs. (\ref{Eq:1D_U012}) and (\ref{Eq:1DLimitedPolynomial}),
a WCC scheme is actually a convex combination of three linear solutions,
so it is linear stable if all the three linear solutions are stable.
As a result, we study the stability of the three linear solutions respectively.
Considering the one-dimensional linear advection equation, i.e., $f=au$ (a is a constant) in Eq. (\ref{Eq:1DHCL}),
then the linear compact central schemes expressed by Eqs. (\ref{Eq:1DIntegral}), (\ref{Eq:Cal_1D_Derivatives}), (\ref{Eq:1DCalU}), and (\ref{Eq:1D_U012})
can be written in the matrix form as
\begin{equation}\label{Eq:1DSchemeMatrixForm}
  \mathbf{q}_{i+1/2}^{n+1/2}=\mathbf{M}_1\cdot\mathbf{q}_{i}^{n}+\mathbf{M}_2\cdot\mathbf{q}_{i+1}^{n},
\end{equation}
where $\mathbf{q}_{i}^{n}=\left[u_i^n,(u_x)_i^n,\cdots,(u_{Px})_i^n\right]^T$ for a $(P+1)$th-order scheme.
Substituting the Fourier series of the solution $\mathbf{q}(x,t^n)=\sum\limits_m\mathbf{A}_m^ne^{ik_mx}$ into Eq. (\ref{Eq:1DSchemeMatrixForm}),
we have
\begin{equation}\label{Eq:1DFourierAnalysis}
  \mathbf{A}_m^{n+1/2}=\left(\mathbf{M}_1e^{-i\theta/2}+\mathbf{M}_2e^{i\theta/2}\right)\cdot\mathbf{A}_m^n=\mathbf{G}\cdot\mathbf{A}_m^n,
\end{equation}
where $\theta=k_m\Delta x\in(0,2\pi]$ is a scaled wave number for convenience and $\mathbf{G}$ is the amplification matrix.
The linear stability condition is
\begin{equation}\label{Eq:LinearStabilityCondition}
  \max\limits_{1\le i\le P}\left|\lambda_i(\mathbf{G})\right|\le 1,
\end{equation}
where $\lambda_i$ is the eigenvalue of $\mathbf{G}$.
The coefficient matrices $\mathbf{M}_1$ and $\mathbf{M}_2$ for 2nd-, 3rd-, and 4th-order compact central schemes are shown in \ref{AppendixA}.
It is difficult to get analytical solutions of the stability conditions because of the complexity of the amplification matrix,
so we use MATLAB to get the numerical solutions as shown by Table \ref{Table:WCC_Linear_Stability}.

\begin{table}
  \centering
  \begin{tabular}{|c|c|c|c|}
     \hline
     &$u_0$   & $u_1$     & $u_2$\\
     \hline
     2nd-order&$\nu\in[0,0.5]$ &$\nu\in[0,0.5]$  &$\nu\in[0,0.5]$  \\
     \hline
     3rd-order&$\nu\in[0,0.384]$ &$\nu\in[0,0.5]$  &$\nu\in[0,0.5]$  \\
     \hline
     4th-order&$\nu\in[0,0.304]$ &$\nu\in[0,0.5]$  &$\nu\in[0,0.5]$  \\
     \hline
  \end{tabular}
  \caption{Linear stability conditions of compact central schemes ($\nu=a\frac{\Delta t}{\Delta x}$).}
  \label{Table:WCC_Linear_Stability}
\end{table}

\section{Arbitrarily high-order WCC schemes in two dimensions}\label{SEC:2DWCCS}
\subsection{Linear compact central schemes in two dimensions}\label{SubSEC:2DLinearWCCS}

We consider the following two-dimensional scalar conservation law, equipped with certain initial conditions and boundary conditions
\begin{subequations}\label{Eq:2DHCL}
  \begin{equation}\label{SubEq:2DHCL_Eq}
    \frac{\partial u}{\partial t}+\frac{\partial f(u)}{\partial x}+\frac{\partial g(u)}{\partial y}=0,\quad (x,y)\in\Omega,\quad t\in[0,\infty),
  \end{equation}
  \begin{equation}\label{SubEq:2DHCL_IC}
    u(x,y,0)=u_0(x,y), \quad (x,y)\in\Omega,
  \end{equation}
  \begin{equation}\label{SubEq:2DHCL_BC}
    u(x,y,t)=u_{\partial\Omega}(x,y,t),\quad (x,y)\in\partial\Omega, \quad t\in[0,\infty).
  \end{equation}
\end{subequations}

\begin{figure}
  \centering
  \includegraphics[width=8 cm]{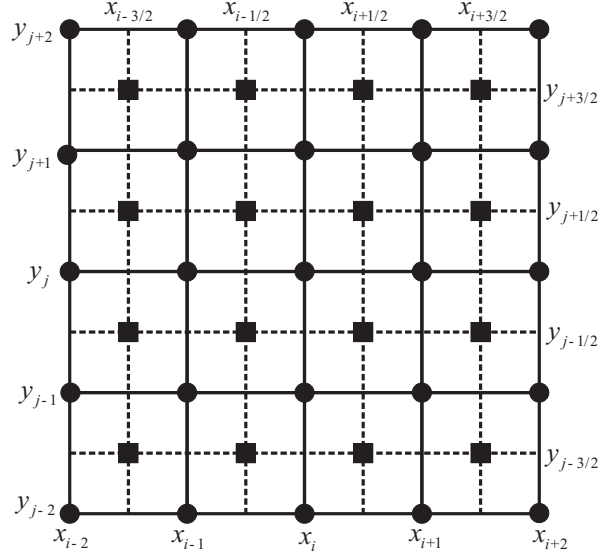}
  \caption{Schematic diagram of the spatial mesh for 2D WCC schemes.
           Solution points at $t^n$ are denoted by $\CIRCLE$ and solution points at $t^{n+1/2}$ are denoted by $\blacksquare$.}
  \label{FIG:2DMesh}
\end{figure}
A regular spatial domain $[x_S,x_E]\times[y_S,y_E]$ is discretized by two sets of Cartesian meshes.
As shown by Fig. \ref{FIG:2DMesh}, the first mesh is defined by the solid lines,
where $x_i=x_S+(i-1)\Delta x$ ($i=1 \text{ to } N_x+1$) with $\Delta x=(x_E-x_S)/N_x$,
and $y_j=y_S+(j-1)\Delta y$ ($j=1 \text{ to } N_y+1$) with $\Delta y=(y_E-y_S)/N_y$.
The other mesh is defined by the dash lines,
where $x_{i+1/2}=(x_i+x_{i+1})/2$ ($i=0 \text{ to } N_x+1$),
$y_{j+1/2}=(y_j+y_{j+1})/2$ ($j=0 \text{ to } N_y+1$).
Similar to the one-dimensional case, a staggered updating strategy is adopted:
\begin{itemize}
  \item At $t^n$, the solution is stored on the original cell $\mathcal{C}_{i,j}$ that is defined by $[x_{i-1/2},x_{i+1/2}]\times[y_{j-1/2},y_{j+1/2}]$.
  \item At $t^{n+1/2}$, the solution is shifted to the staggered cell $\mathcal{C}_{i+1/2,j+1/2}$ that is defined by $[x_{i},x_{i+1}]\times[y_{j},y_{j+1}]$.
  \item At $t^{n+1}$, the solution reverts to $\mathcal{C}_{i,j}$, entering the next loop.
\end{itemize}
On cell $\mathcal{C}_{i,j}$ at $t^{n}$, the solution is approximated by a $P$th-order polynomial as
\begin{equation}\label{Eq:2DPolynomial_Space}
  u(x,y,t^{n})=\sum_{l=0}^{P}\sum_{k=0}^{P-l}\frac{1}{k!}\frac{1}{l!}(u_{kx,ly})_{i,j}^{n}\left(\frac{x-x_i}{\Delta x}\right)^k\left(\frac{y-y_j}{\Delta y}\right)^l,
  \quad (x,y)\in \mathcal{C}_{i,j},
\end{equation}
where $(u_{kx,ly})_{i,j}^{n}$ is the scaled spatial derivatives of the approximated solution at the cell center $(x_i,y_j,t^n)$.

\begin{figure}
  \centering
  \includegraphics[width=8 cm]{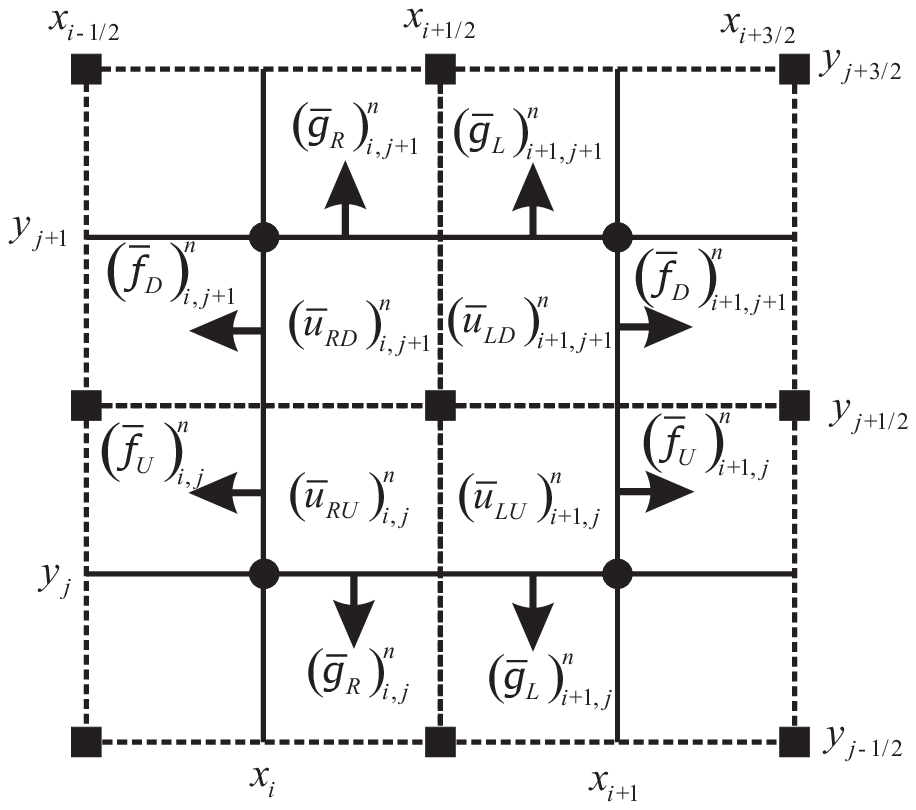}
  \caption{Schematic diagram of 2D fluxes in Eq. (\ref{Eq:2DIntegral}).
           Solution points at $t^n$ are denoted by $\CIRCLE$ and solution points at $t^{n+1/2}$ are denoted by $\blacksquare$.}
  \label{FIG:2DFluxes}
\end{figure}

The time marching scheme from $t=t^n$ to $t=t^{n+1/2}$ is almost the same as that from $t=t^{n+1/2}$ to $t=t^{n+1}$,
except that the mesh is shifted by $\left(\Delta x/2,\Delta y/2\right)$.
Without loss of generality, we consider the time marching scheme from $t=t^n$ to $t=t^{n+1/2}$.
First, we perform the integral of Eq. (\ref{Eq:2DHCL}) on the control volume $[x_{i},x_{i+1}]\times[y_{j},y_{j+1}]\times[t^n,t^{n+1/2}]$,
obtaining the cell average by a classical finite-volume scheme as
\begin{equation}\label{Eq:2DIntegral}
  \begin{split}
    \bar{u}_{i+1/2,j+1/2}^{n+1/2}&=\frac{1}{4}\left[(\bar{u}_{RU})_{i,j}^n+(\bar{u}_{RD})_{i,j+1}^n+(\bar{u}_{LU})_{i+1,j}^n+(\bar{u}_{LD})_{i+1,j+1}^n\right]\\
    &+\frac{\Delta t}{2\Delta x}\left[(\bar{f}_U)_{i,j}^n+(\bar{f}_D)_{i,j+1}^n-(\bar{f}_U)_{i+1,j}^n-(\bar{f}_D)_{i+1,j+1}^n\right]\\
    &+\frac{\Delta t}{2\Delta y}\left[(\bar{g}_R)_{i,j}^n+(\bar{g}_L)_{i+1,j}^n-(\bar{g}_R)_{i,j+1}^n-(\bar{g}_L)_{i+1,j+1}^n\right],
  \end{split}
  \end{equation}
where all the terms on the right hand side represents the fluxes through the control volume.
The physical meaning of each term is intuitively shown in Fig \ref{FIG:2DFluxes},
and the mathematical form of each term is given in \ref{AppendixB}.

We note that $f$, $g$ and their derivatives can be calculated by using $u$ and the derivatives of $u$,
and the time derivatives can be calculated by repeatedly utilizing the original PDE, Eq. (\ref{SubEq:2DHCL_Eq}).
Therefore, only $u$ and its spatial derivatives are independent variables.
Next, we use a central difference to calculate the spatial derivatives.
The Taylor series at the four vertices of $\mathcal{C}_{i+1/2,j+1/2}$ read
\begin{subequations}\label{Eq:2DTaylorSeries}
  \begin{equation}\label{SubEq:2DTaylorSeries1}
   (u_{kx,ly})_{i,j}^{n+1/2}=\sum_{t=0}^{P-k-l}\sum_{s=0}^{P-k-l-t}\frac{1}{s!}\frac{1}{t!}(u_{(k+s)x,(l+t)y})_{i+1/2,j+1/2}^{n+1/2}\left(-\frac{1}{2}\right)^{s+t},
  \end{equation}

  \begin{equation}\label{SubEq:2DTaylorSeries2}
   (u_{kx,ly})_{i+1,j}^{n+1/2}=\sum_{t=0}^{P-k-l}\sum_{s=0}^{P-k-l-t}\frac{1}{s!}\frac{1}{t!}(u_{(k+s)x,(l+t)y})_{i+1/2,j+1/2}^{n+1/2}\frac{(-1)^t}{2^{s+t}},
  \end{equation}

  \begin{equation}\label{SubEq:2DTaylorSeries3}
   (u_{kx,ly})_{i+1,j+1}^{n+1/2}=\sum_{t=0}^{P-k-l}\sum_{s=0}^{P-k-l-t}\frac{1}{s!}\frac{1}{t!}(u_{(k+s)x,(l+t)y})_{i+1/2,j+1/2}^{n+1/2}\frac{1}{2^{s+t}},
  \end{equation}

  \begin{equation}\label{SubEq:2DTaylorSeries4}
   (u_{kx,ly})_{i,j+1}^{n+1/2}=\sum_{t=0}^{P-k-l}\sum_{s=0}^{P-k-l-t}\frac{1}{s!}\frac{1}{t!}(u_{(k+s)x,(l+t)y})_{i+1/2,j+1/2}^{n+1/2}\frac{(-1)^s}{2^{s+t}}.
  \end{equation}
\end{subequations}
(\ref{SubEq:2DTaylorSeries2}+\ref{SubEq:2DTaylorSeries3})-(\ref{SubEq:2DTaylorSeries1}+\ref{SubEq:2DTaylorSeries4}) yields
\begin{subequations}\label{Eq:2DDerivatives}
  \begin{equation}\label{SubEq:2DUx}
  \begin{split}
     &(u_{(k+1)x,ly})_{i+1/2,j+1/2}^{n+1/2}=0.5\left[(u_{kx,ly})_{i+1,j}^{n+1/2}+(u_{kx,ly})_{i+1,j+1}^{n+1/2}-(u_{kx,ly})_{i,j}^{n+1/2}\right.\\
     &\left.-(u_{kx,ly})_{i,j+1}^{n+1/2}\right]-\sum_{s=3}^{P-k-l}\frac{1}{s!}(u_{(k+s)x,ly})_{i+1/2,j+1/2}^{n+1/2}\left[\frac{1}{2^{s}}-\left(-\frac{1}{2}\right)^s\right]\\
     &-0.5\sum_{t=2}^{P-k-l}\sum_{s=1}^{P-k-l-t}\frac{1}{s!}\frac{1}{t!}(u_{(k+s)x,(l+t)y})_{i+1/2,j+1/2}^{n+1/2}\frac{\left[1-(-1)^s\right]\left[1+(-1)^t\right]}{2^{s+t}},
  \end{split}
  \end{equation}
  \text{(\ref{SubEq:2DTaylorSeries3}+\ref{SubEq:2DTaylorSeries4})-(\ref{SubEq:2DTaylorSeries1}+\ref{SubEq:2DTaylorSeries2}) yields}
  \begin{equation}\label{SubEq:2DUy}
  \begin{split}
     &(u_{kx,(l+1)y})_{i+1/2,j+1/2}^{n+1/2}=0.5\left[(u_{kx,ly})_{i,j+1}^{n+1/2}+(u_{kx,ly})_{i+1,j+1}^{n+1/2}-(u_{kx,ly})_{i,j}^{n+1/2}\right.\\
     &\left.-(u_{kx,ly})_{i+1,j}^{n+1/2}\right]-\sum_{t=3}^{P-k-l}\frac{1}{t!}(u_{kx,(l+t)y})_{i+1/2,j+1/2}^{n+1/2}\left[\frac{1}{2^{t}}-\left(-\frac{1}{2}\right)^t\right]\\
     &-0.5\sum_{t=1}^{P-k-l}\sum_{s=2}^{P-k-l-t}\frac{1}{s!}\frac{1}{t!}(u_{(k+s)x,(l+t)y})_{i+1/2,j+1/2}^{n+1/2}\frac{\left[1+(-1)^s\right]\left[1-(-1)^t\right]}{2^{s+t}},
  \end{split}
  \end{equation}
  \text{where}
  \begin{equation}\label{SubEq:2DUij}
    (u_{kx,ly})_{i,j}^{n+1/2}=\sum_{m=0}^{P-k-l} \frac{1}{m!}(u_{kx,ly,mt})_{i,j}^{n}
  \end{equation}
\end{subequations}
When $P-k-l<3$ the last two summation terms of Eqs. (\ref{SubEq:2DUx}) and (\ref{SubEq:2DUy}) vanish,
otherwise they contain high-order derivatives at $t^{n+1/2}$.
Therefore, we have to update the derivatives from high-order to low-order to obtain an explicit scheme.

\begin{remark}
  We may get different numerical solutions for cross-derivative terms by using Eqs. (\ref{SubEq:2DUx}) and (\ref{SubEq:2DUy}), e.g., $u_{xy}\ne u_{yx}$.
  In order to save computing memories, we take the arithmetic mean as the final solution,
  e.g., $u_{x,y}=0.5(u_{xy}+u_{yx})$.
\end{remark}

\subsection{Weighted compact central schemes in two dimensions}\label{SubSEC:2DWCC}
Similar to the one-dimensional case, when a solution contains discontinuities,
we design a WENO-type limiter for two-dimensional compact central schemes as
\begin{equation}\label{Eq:2DLimitedPolynomial}
  u_*^n(x,y)=\sum_{m=0}^{4}w_mu_m^n(x,y) \text{ with }w_m=\frac{\hat{w}_m}{\sum_{k=0}^{4}\hat{w}_k}, \quad x\in \mathcal{C}_{i,j},
\end{equation}
where $u_0(x,y)$ denotes the high-order solution derived by 2D linear compact central schemes
\begin{equation}\label{Eq:2DHighOrderPolynomial}
  u_0^n(x,y)=\sum_{l=0}^{P}\sum_{k=0}^{P-l}\frac{1}{k!}\frac{1}{l!}(u_{kx,ly}^0)_{i,j}^{n}\left(\frac{x-x_i}{\Delta x}\right)^k\left(\frac{y-y_j}{\Delta y}\right)^l,
\end{equation}
and $u_m^n(x,y)$ ($m=1\text{ to } 4$) are four additional first-order polynomial constructed based
on the four edges of $\mathcal{C}_{i,j}$.
The specific forms of $u_m^n(x,y)$ ($m=1\text{ to } 4$) are respectively constructed as
\begin{equation}\label{Eq:2DLowOrderPolynomial}
  u_m^n(x,y)=(u^m)_{i,j}^{n}+(u_x^m)_{i,j}^{n}\left(\frac{x-x_i}{\Delta x}\right)+(u_y^m)_{i,j}^{n}\left(\frac{y-y_j}{\Delta y}\right),
\end{equation}
where
\begin{equation*}
  (u^1)_{i,j}^{n}=(u^2)_{i,j}^{n}=(u^3)_{i,j}^{n}=(u^4)_{i,j}^{n}=\overline{u_0^n(x,y)},
\end{equation*}
\begin{equation*}
  (u_x^1)_{i,j}^{n}=2(u^1)_{i,j}^{n}-(u_{i-1/2,j-1/2}^n+u_{i-1/2,j+1/2}^n),
\end{equation*}
\begin{equation*}
  (u_y^1)_{i,j}^{n}=u_{i-1/2,j+1/2}^n-u_{i-1/2,j-1/2}^n,
\end{equation*}
\begin{equation*}
  (u_x^2)_{i,j}^{n}=u_{i+1/2,j-1/2}^n-u_{i-1/2,j-1/2}^n,
\end{equation*}
\begin{equation*}
  (u_y^2)_{i,j}^{n}=2(u^2)_{i,j}^{n}-(u_{i-1/2,j-1/2}^n+u_{i+1/2,j-1/2}^n),
\end{equation*}
\begin{equation*}
  (u_x^3)_{i,j}^{n}=(u_{i+1/2,j-1/2}^n+u_{i+1/2,j+1/2}^n)-2(u^3)_{i,j}^{n},
\end{equation*}
\begin{equation*}
  (u_y^3)_{i,j}^{n}=u_{i+1/2,j+1/2}^n-u_{i+1/2,j-1/2}^n,
\end{equation*}
\begin{equation*}
  (u_x^4)_{i,j}^{n}=u_{i+1/2,j+1/2}^n-u_{i-1/2,j+1/2}^n,
\end{equation*}
\begin{equation*}
  (u_y^4)_{i,j}^{n}=(u_{i-1/2,j+1/2}^n+u_{i+1/2,j+1/2}^n)-2(u^4)_{i,j}^{n}.
\end{equation*}
We note that $\overline{u_0^n(x,y)}$ is derived by Eq. (\ref{Eq:2DIntegral}),
and $u_{i\pm1/2,j\pm1/2}^n$ at the four vertices of $\mathcal{C}_{i,j}$ are derived by Eq. (\ref{SubEq:2DUij}),
so the construction of $u_m^n(x,y)$ ($m=1\text{ to } 4$) is straightforward.

Depending on the order of the scheme, the weights are calculated in the following two ways:
\begin{enumerate}
  \item In the case of the second-order scheme ($P=1$), $u_m^n(x,y)$ $(m=0\text{ to }4)$ have a uniform degree.
  The weights are calculated as
  \begin{equation}\label{Eq:2D2ndOrderWeights}
    \hat{w}_m=\left(\frac{1}{\beta_m+\epsilon}\right)^\alpha,
  \end{equation}
  where $\epsilon=10^{-40}$ is used to avoid dividing by 0,
  and $\beta_m$ is the smoothness indicator of $u_m^n(x,y)$.
  Extending from the 1D case, the smoothness indicator of a 2D $P$th-order polynomial defined on $\mathcal{C}_{i,j}$ is calculated by
  \begin{equation}\label{Eq:2DSmoothIndicator}
    \beta_m=\sum_{l=0}^{P}\sum_{k=0,k+l>0}^{P-l}\Delta x^{2k-1}\Delta y^{2l-1}\int_{x_i-\Delta x/2}^{x_i+\Delta x/2}
    \int_{y_j-\Delta y/2}^{y_j+\Delta y/2} \left(\frac{\partial^{k+l}u_m^n(x,y)}{\partial x^k\partial y^l}\right)^2dxdy.
  \end{equation}

  \item In the case of the higher than second-order schemes ($P\ge2$), $u_m^n(x,y)$ $(m=0\text{ to }4)$ have unequal degrees.
  The weights are calculated as
  \begin{subequations}\label{Eq:1DHighOrderWeights}
    \begin{equation}
    \hat{w}_0=1+\left(\frac{\sigma\tau}{\beta_0^2+\epsilon}\right)^\alpha,\quad \hat{w}_m=\left(\frac{\sigma\tau}{\beta_m^2+\epsilon}\right)^\alpha (m=1 \text{ to }4),
   \end{equation}
   \text{where}
    \begin{equation}
    \begin{split}
       \sigma&=\left[u_0^n\left(x_{i-1/2},y_{j-1/2}\right)-u_{i-1/2,j-1/2}^n\right]^2\\
             &+\left[u_0^n\left(x_{i+1/2},y_{j-1/2}\right)-u_{i+1/2,j-1/2}^n\right]^2\\
             &+\left[u_0^n\left(x_{i+1/2},y_{j+1/2}\right)-u_{i+1/2,j+1/2}^n\right]^2\\
             &+\left[u_0^n\left(x_{i-1/2},y_{j+1/2}\right)-u_{i-1/2,j+1/2}^n\right]^2,
    \end{split}
   \end{equation}
   \begin{equation}
   \begin{split}
      \tau&=\sum_{k=0}^{P}\Delta x^{2k-1}\Delta y^{2(P-k)-1}\int_{x_i-\Delta x/2}^{x_i+\Delta x/2}\int_{y_j-\Delta y/2}^{y_j+\Delta y/2}
    \left(\frac{d^Pu_0^n(x,y)}{\partial x^k\partial y^{P-k}}\right)^2dxdy\\
     &=\sum_{k=0}^{P}\left[(u_{kx,(P-k)y}^0)_{i,j}^n\right]^2,
   \end{split}
   \end{equation}
  \end{subequations}
   where, $u_0^n\left(x_{i\pm1/2},y_{j\pm1/2}\right)$ is derived by a space approximation using Eq. (\ref{Eq:2DHighOrderPolynomial}),
   while $u_{i\pm1/2,j\pm1/2}^n$ is derived by a time approximation using Eq. (\ref{SubEq:2DUij}).
   Therefore, $\sigma$ is the measurement of the residuals at the four vertices of $\mathcal{C}_{i,j}$,
   and $\tau$ is the smoothness indicator for the highest-order terms of $u_0^n(x,y)$.
\end{enumerate}

As wee can see, the 2D limiter is a straightforward extension of the 1D limiter,
so the order of convergence at smooth regions and shock-capturing properties follow exactly the same as the one-dimensional case.

\section{High-Order WCC schemes for hyperbolic systems}\label{SEC:_System_Extension}
In Sec. \ref{SEC:1DWCCS} and Sec. \ref{SEC:2DWCCS},
we introduce WCC schemes for solving scalar hyperbolic conservation laws in one dimension and two dimensions respectively.
When we extend the WCC schemes to solve hyperbolic systems,
almost every thing follows exactly the same as the scalar cases,
except the Cauchy–Kovalewski procedure which seems to become cumbersome for high-order cases.
Fortunately, there is a very efficient algorithm based on the generalized Leibniz rule
which was originally developed by Dyson \cite{Dyson2002Cauchy_Kovalewski},
and was further developed by Dumbser and Munz for ADER schemes \cite{Dumbser2006Cauchy_Kovalewski}.
Although, this technique is originally developed for the Euler equations,
it can be extended to other systems using similar technique as well.

In addition, when we implement the WENO-type limiter to a hyperbolic system,
we can choose different variables to apply the limiter.
Usually, the most straightforward and efficient choice is conservative variables,
but the numerical oscillations near strong discontinuities can not be completely eliminated by high-order WCC schemes.
Following many other high-order schemes, we apply the limiter to characteristic variables to effectively control spurious oscillations.

For a one-dimensional system
\begin{equation}\label{Eq:1DHCLSystem}
  \frac{\partial \mathbf{U}}{\partial t}+\frac{\partial \mathbf{F}}{\partial x}=0,
\end{equation}
the limiting process is as follows:
\begin{itemize}
  \item Calculate the left eigenvector matrix $\mathbf{L(U)}$ and the right eigenvector matrix $\mathbf{R(U)}$
  of the Jacobian matrix $\mathbf{A(U)}=\frac{\partial \mathbf{F}}{\partial \mathbf{U}}$, i.e. $\mathbf{A(U)}=\mathbf{R(U)}\mathbf{\Lambda}\mathbf{L(U)}$.
  Here, $\mathbf{\Lambda}$ is a diagonal matrix with $\mathbf{\Lambda}_{kk}=\lambda_k$, where $\lambda_k$ is the eigenvalue of $\mathbf{A(U)}$;
  \item Project the three solutions in Eq. (\ref{Eq:1D_U012}) into the characteristic space as
  $\tilde{\mathbf{U}}_m^n(x)=\mathbf{L}(\mathbf{\bar{U}}_i^n)\cdot\mathbf{U}_m^n(x)$;
  \item Apply the limiter to the characteristic vector $\tilde{\mathbf{U}}_m^n(x)$ component by component as
  \begin{equation*}
  (\tilde{\mathbf{U}}_k)_*^n(x)=\sum_{m=0}^{2}\tilde{w}_m^k(\tilde{\mathbf{U}}_k)_m^n(x);
\end{equation*}
  \item Project the limited characteristic vector $\tilde{\mathbf{U}}_*^n(x)$ back to physical space as
  $\mathbf{U}_*^n(x)=\mathbf{R}(\mathbf{\bar{U}}_i^n)\cdot\tilde{\mathbf{U}}_*^n(x)$.
\end{itemize}

For a two-dimensional system
\begin{equation}\label{Eq:2DHCLSystem}
  \frac{\partial \mathbf{U}}{\partial t}+\frac{\partial \mathbf{F}}{\partial x}+\frac{\partial \mathbf{G}}{\partial y}=0,
\end{equation}
there are two normal directions in each cell of a Cartesian mesh.
A traditional method performs characteristic decompositions in every normal direction
and applies the limiter two times.
This process consumes a lot of CPU times and makes the algorithm inefficient.
Here, we adopt a rotated characteristic decomposition technique \cite{Shen2020RotatedDecompositon} as follows:
\begin{itemize}
  \item Determine the rotated direction as
  \begin{equation}\label{Eq:RotatedNormal}
    \mathbf{n}=\frac{\nabla\phi}{|\nabla\phi|}=(cos\theta,sin\theta)
  \end{equation}
  where $\nabla\phi$ is the gradient of a physical quantity $\phi$ which is approximated by a central difference;

  \item Calculate the left eigenvector matrix $\mathbf{L}(\mathbf{U},\theta)$ and the right eigenvector matrix $\mathbf{R}(\mathbf{U},\theta)$
  of the Jacobian matrix $\mathbf{A}(\mathbf{U},\theta)=\frac{\partial \mathbf{F}}{\partial \mathbf{U}}cos\theta+\frac{\partial \mathbf{G}}{\partial \mathbf{U}}sin\theta$;

  \item Project the linear solutions in Eqs. (\ref{Eq:2DHighOrderPolynomial}) and (\ref{Eq:2DLowOrderPolynomial}) into the characteristic space as
  $\tilde{\mathbf{U}}_m^n(x,y)=\mathbf{L}(\mathbf{\bar{U}}_{i,j}^n,\theta_{i,j})\cdot\mathbf{U}_m^n(x,y)$;

  \item Apply the limiter to the characteristic vector $\tilde{\mathbf{U}}_m^n(x,y)$ component by component as
  \begin{equation*}
  (\tilde{\mathbf{U}}_k)_*^n(x,y)=\sum_{m=0}^{4}\tilde{w}_m^k(\tilde{\mathbf{U}}_k)_m^n(x,y);
\end{equation*}
  \item Project the limited characteristic vector $\tilde{\mathbf{U}}_*^n(x,y)$ back to physical space as
  $\mathbf{U}_*^n(x,y)=\mathbf{R}(\mathbf{\bar{U}}_{i,j}^n,\theta_{i,j})\cdot\tilde{\mathbf{U}}_*^n(x,y)$.
\end{itemize}

As we can see, this technique requires only one-time characteristic decomposition in the direction of a physical quantity's gradient,
so it reduces the computational cost significantly.
When we implemented a third-order finite-volume WENO scheme for solving the Euler equations on Cartesian meshes,
the rotated characteristic decomposition technique eliminates spurious oscillations effectively
and saves about 40\% CPU times comparing to the traditional characteristic decomposition method \cite{Shen2020RotatedDecompositon}.
If we want to further reduce the computational cost,
we can employ the current approach in addition to the existing ones
that only perform characteristic decomposition near discontinuities
\cite{Ren2003hybridCompatWENO, Puppo2003adaptive, Puppo2011adaptive, Li2010hybridWENO, Peng2019adaptiveWENOZ}.

\section{Numerical examples}\label{SEC:_NumExam}
For convenience, we denote $P$th-order linear and weighted compact central schemes as LCCS-$P$ and WCCS-$P$ respectively.
In the light of the linear stability analysis in Sec. \ref{SubSEC:1DLinearStability},
we set CFL=0.4 for LCCS-2 and WCCS-2, CFL=0.3 for LCCS-3 and WCCS-3, CFL=0.25 for LCCS-4 and WCCS-4 in all the numerical examples.

\subsection{Numerical examples for the 1D linear advection equation}
We consider the 1D linear advection equation
\begin{equation}\label{Eq:1DAdvEq}
  \frac{\partial u}{\partial t}+\frac{\partial u}{\partial x}=0.
\end{equation}
The first case is the advection of a sinusoidal wave. The computational domain is $[-1,1]$,
the initial condition is $u_0(x)=sin(\pi x)$,
and periodic boundary conditions are applied on the left and right sides.
Table \ref{Table:1DAdvectionConvergence} shows the numerical convergence of LCCS and WCCS at $t=2$,
where $L_1$ is the average error and $L_\infty$ is the maximum error.
We observe that LCCS achieve the expected convergence rate for $L_1$ and $L_\infty$.
When applied the limiter, WCCS-2 loses some accuracy for $L_\infty$
because the TVD-type limiter degrades the second-order scheme to first-order at local extrema.
However, the WENO-type limiter has no effect on the 3rd- and 4th-order schemes.
The numerical results agree well with the analysis in Sec. \ref{SubSEC:1DWCC}.

\begin{table}
  \centering
  \begin{tabular}{|c|c|c|c|c|c|c|c|c|}
    \hline
      & \multicolumn{4}{|c|}{LCCS-2} & \multicolumn{4}{|c|}{WCCS-2} \\
     \hline
      Mesh size       & $L_1$  &Order & $L_\infty$  &Order    & $L_1$  &Order & $L_\infty$  &Order    \\
     \hline
     1/25              &9.66E-4 &-     &1.50E-3      &-        &2.18E-3 &-     &9.54E-3      &- \\

     1/50             &2.38E-4 &2.02  &3.72E-4      &2.01     &6.10E-4 &1.84  &3.80E-3      &1.33 \\

     1/100             &5.94E-5 &2.00  &9.30E-5      &2.00     &1.54E-4 &1.99  &1.71E-3      &1.15 \\

     1/200             &1.48E-5 &2.00  &2.33E-5      &2.00     &3.80E-5 &2.02  &6.43E-4      &1.41 \\
     \hline

     & \multicolumn{4}{|c|}{LCCS-3} & \multicolumn{4}{|c|}{WCCS-3} \\
     \hline
      Mesh size       & $L_1$  &Order & $L_\infty$  &Order    & $L_1$  &Order & $L_\infty$  &Order    \\
     \hline
     1/25              &2.77E-5 &-     &4.42E-5      &-        &2.77E-5 &-     &4.42E-5      &-  \\

     1/50             &3.49E-6 &2.99  &5.53E-6      &3.00     &3.49E-6 &2.99  &5.53E-6      &3.00 \\

     1/100             &4.39E-7 &2.99  &6.92E-7      &3.00     &4.39E-7 &2.99  &6.92E-7      &3.00 \\

     1/200             &5.49E-8 &3.00  &8.65E-8      &3.00     &5.49E-8 &3.00  &8.65E-8      &3.00 \\
     \hline

     & \multicolumn{4}{|c|}{LCCS-4} & \multicolumn{4}{|c|}{WCCS-4} \\
     \hline
      Mesh size       & $L_1$  &Order & $L_\infty$  &Order    & $L_1$  &Order & $L_\infty$  &Order    \\
     \hline
     1/25              &8.75E-7 &-     &1.36E-6      &-        &8.75E-7 &-     &1.36E-6      &-  \\

     1/50             &5.43E-8 &4.01  &8.49E-8      &4.00     &5.43E-8 &4.01  &8.49E-8      &4.00 \\

     1/100             &3.37E-9 &4.01  &5.28E-9      &4.01     &3.37E-9 &4.01  &5.28E-9      &4.01 \\

     1/200             &2.10E-10 &4.00  &3.29E-10      &4.00     &2.10E-10 &4.00  &3.29E-10      &4.00  \\
     \hline
  \end{tabular}
  \caption{Numerical errors of the advection of a sinusoidal wave at $t=2$ computed by LCCS and WCCS with different mesh sizes.}
  \label{Table:1DAdvectionConvergence}
\end{table}

\begin{figure}
  \centering
  \subfigure[WCCS-2]{
  \label{FIG:Linear_Advection_WCCS2}
  \includegraphics[width=7 cm]{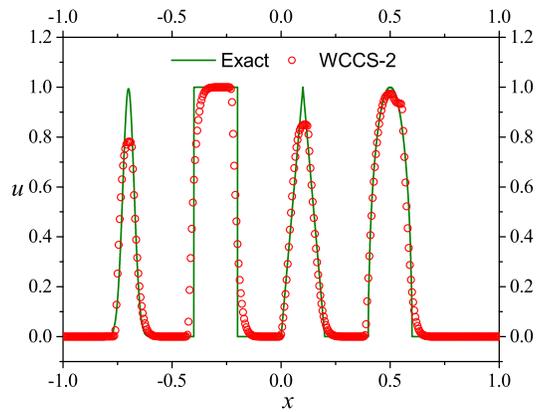}}
  \subfigure[WCCS-3]{
  \label{FIG:Linear_Advection_WCCS3}
  \includegraphics[width=7 cm]{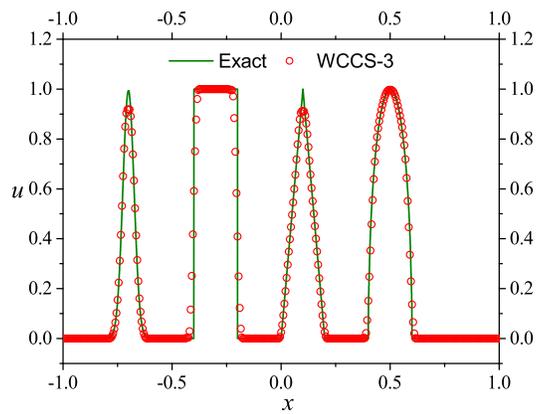}}
  \subfigure[WCCS-4]{
  \label{FIG:Linear_Advection_WCCS4}
  \includegraphics[width=7 cm]{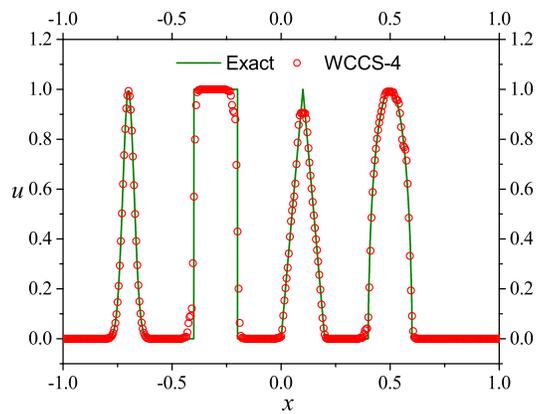}}
  \caption{The advection of a combination of Gaussians, a square wave, a sharp triangle wave,
  and a half ellipse arranged from left to right at $t=12$ calculated by WCCS with $\Delta x=1/200$.}
\label{FIG:Linear_Advection}
\end{figure}

The second case is the advection of a combination of Gaussians, a square wave, a sharp
triangle wave, and a half ellipse arranged from left to right.
This case was originally proposed by Jiang and Shu \cite{Jiang_Shu1996WENO},
and then was widely used to test the high-fidelity property high-order schemes for different shape of solutions.
The computational domain is $[-1,1]$.
The initial condition is given by
\begin{subequations}
\begin{equation*}
  u_0(x)=\begin{cases}
        \frac{1}{6}\left[G(x,\beta,z-\delta)+4G(x,\beta,z)+G(x,\beta,z+\delta)\right], & \mbox{if } -0.8\le x\le -0.6 \\
        1, & \mbox{if } -0.4\le x\le -0.2 \\
        1-10\left|x-0.1\right|, & \mbox{if } 0\le x\le 0.2 \\
        \frac{1}{6}\left[F(x,\alpha,a-\delta)+4F(x,\alpha,a)+F(x,\alpha,a+\delta)\right], & \mbox{if } 0.4\le x\le 0.6 \\
        0, & \mbox{otherwise}.
  \end{cases}
\end{equation*}
\text{where,}\\
\begin{equation*}
  G(x,\beta,z)=e^{-\beta(x-z)^2},
\end{equation*}

\begin{equation*}
  F(x,\alpha,a)=\sqrt{\mbox{max}(1-\alpha^2(x-a)^2,0)}.
\end{equation*}

\end{subequations}
The constants are set as $a=0.5$, $z=-0.7$, $\delta=0.005$, $\alpha=10$, and $\beta=\frac{\log2}{36\delta^2}$.
Periodic boundary conditions are implemented on the left and right sides.
We run this case to $t=12$ with 400 uniform cells.
Fig. \ref{FIG:Linear_Advection} shows the profiles of $u$ calculated by WCCS comparing with the exact solution.
We observe that WCCS-3 and WCCS-4 perform much better for all kinds of solutions than WCCS-2.
The solution of WCCS-2 has an evident kink for the half ellipse which can be eliminated by refining the mesh.
WCCS-4 has the best performance for the smooth Gaussians among the three schemes,
but WCCS-3 has a shaper transition for discontinuous solutions than WCCS-4.
This is because the smoothness indicator Eq. (\ref{Eq:1DSmoothIndicator}) of a 3rd-order polynomial
is much larger than that of a 4th-order polynomial at discontinuities.
As a result, the limited polynomial Eq. (\ref{Eq:1DLimitedPolynomial}) of WCCS-4
is closer to the 1st-order polynomials than that of WCCS-3.

\subsection{Numerical examples for the 1D Euler equations}
The 1D compressible Euler equations can be written as Eq. (\ref{Eq:1DHCLSystem})
with $\mathbf{U}=[\rho,\rho u,\rho e]^T$ and $\mathbf{F}=[{\rho u,\rho u^2+p,(\rho e+p)u}]^T$,
where $\rho$, $u$, $p$, and $e$ denote the density, velocity, pressure, and specific total energy respectively.
We use the prefect gas model and the specific total energy is calculated as $e=\frac{p}{\rho(\gamma-1)}+\frac{1}{2} u^2$.
Here, we set the specific heat ratio $\gamma=1.4$.

\begin{figure}
  \centering
  \subfigure[WCCS-2]{
  \label{FIG:Sod_WCCS2}
  \includegraphics[width=7 cm]{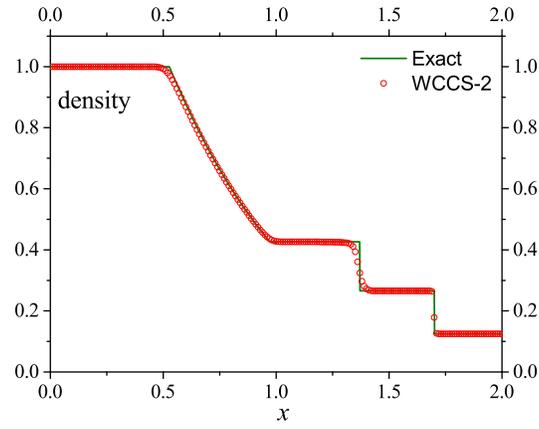}}
  \subfigure[WCCS-3]{
  \label{FIG:Sod_WCCS3}
  \includegraphics[width=7 cm]{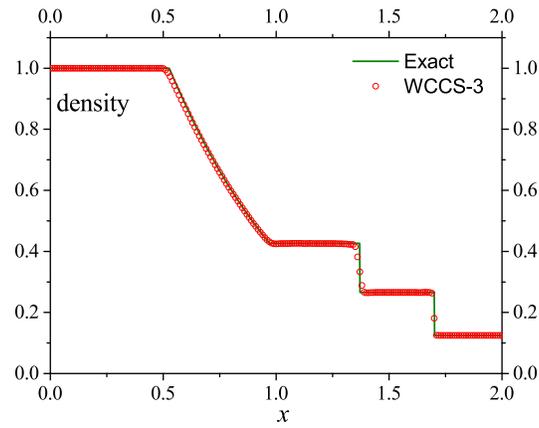}}
  \subfigure[WCCS-4]{
  \label{FIG:Sod_WCCS4}
  \includegraphics[width=7 cm]{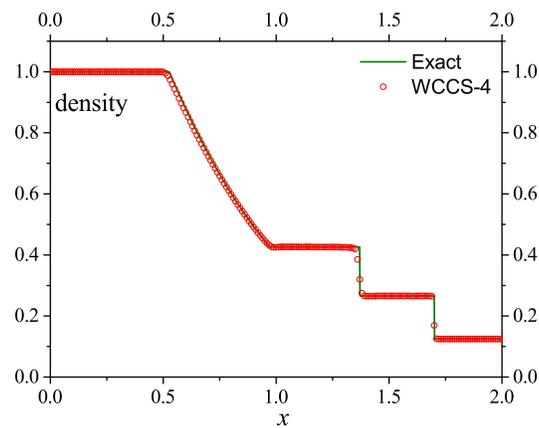}}
  \caption{The density profiles of the Sod shock tube problem at $t=0.4$ calculated by WCCS with $\Delta x=1/100$.}
\label{FIG:Sod}
\end{figure}

The first case is the Sod \cite{Sod1978Survey} shock tube problem that is used to test the shock-capturing property of numerical methods.
The computational domain is [0,2], and the initial condition is given by
\begin{equation*}
  (\rho,u,p)=\begin{cases}
               (1,0,1), & \mbox{if } x<1 \\
               (0.125,0,0.1), & \mbox{otherwise}.
             \end{cases}
\end{equation*}
Non-reflection boundary conditions are applied on the left and right sides.
Fig. \ref{FIG:Sod} shows the density profiles at $t=0.4$ calculated by WCCS with $\Delta x=1/100$.
As we see, WCCS-3 and WCCS-4 obviously perform better than WCCS-2, especially at the contact discontinuity.
WCCS-3 and WCCS-4 have similar performance, because the simple structures do not fully exploit the advantage of high-order schemes.

\begin{figure}
  \centering
  \subfigure[WCCS-2]{
  \label{FIG:Titarev_Toro_WCCS2}
  \includegraphics[width=8 cm]{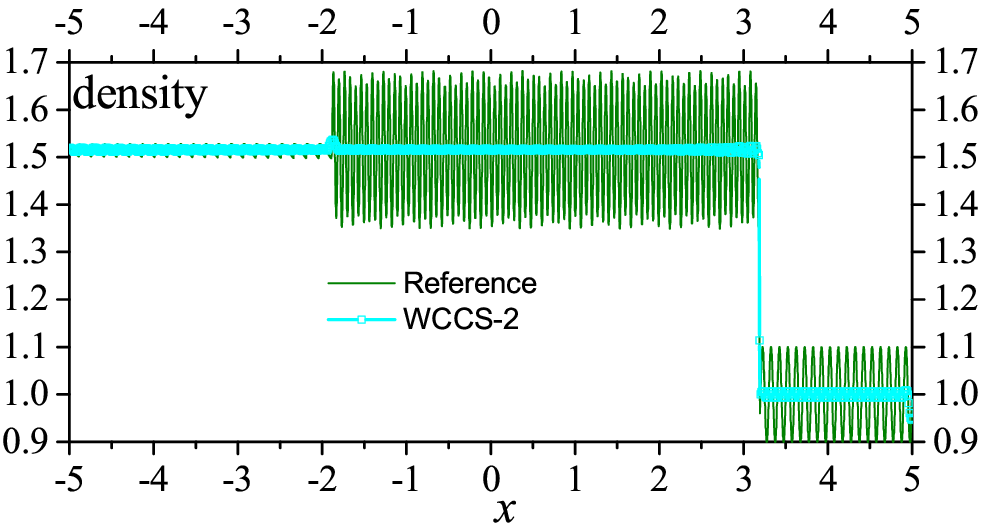}}
  \subfigure[WCCS-3]{
  \label{FIG:Titarev_Toro_WCCS3}
  \includegraphics[width=8 cm]{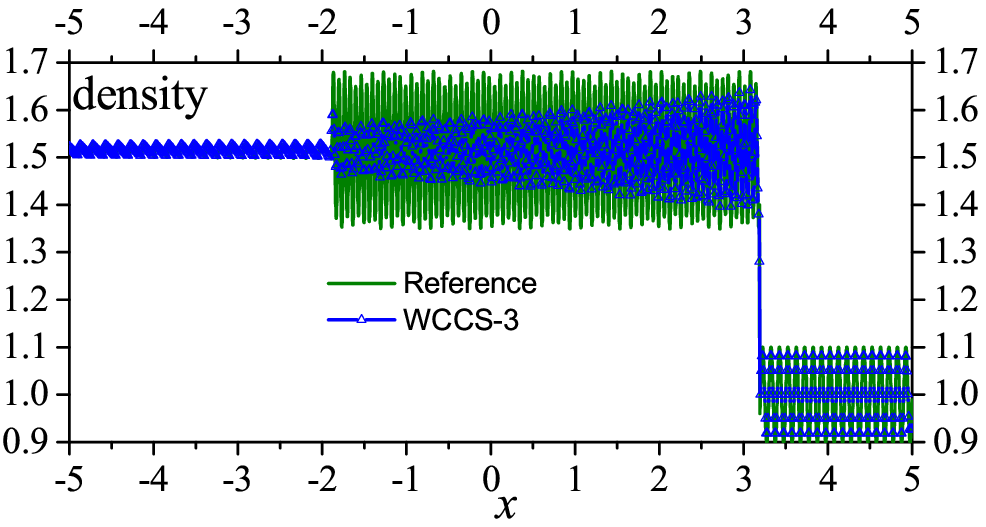}}
  \subfigure[WCCS-4]{
  \label{FIG:Titarev_Toro_WCCS4}
  \includegraphics[width=8 cm]{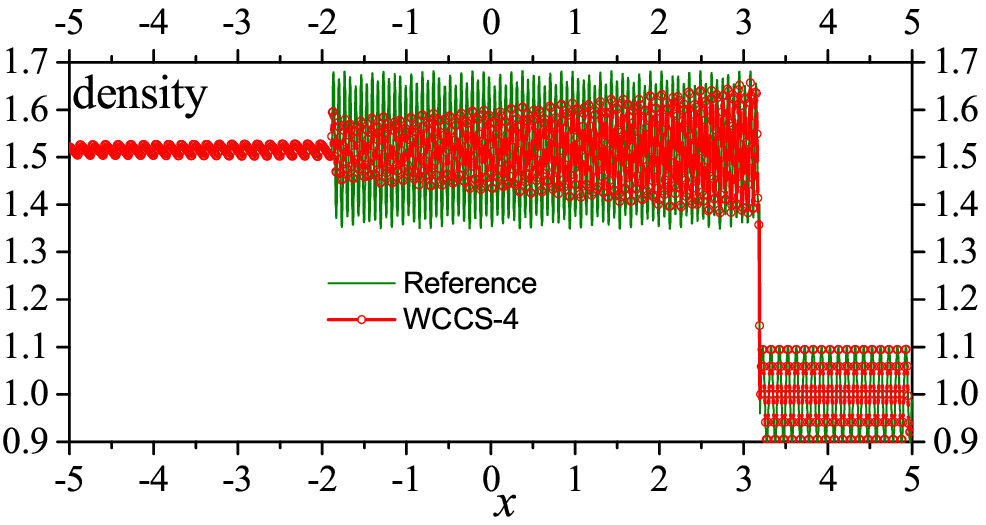}}
  \subfigure[Comparison of the enlarged view]{
  \label{FIG:Titarev_Toro_Comparison}
  \includegraphics[width=8 cm]{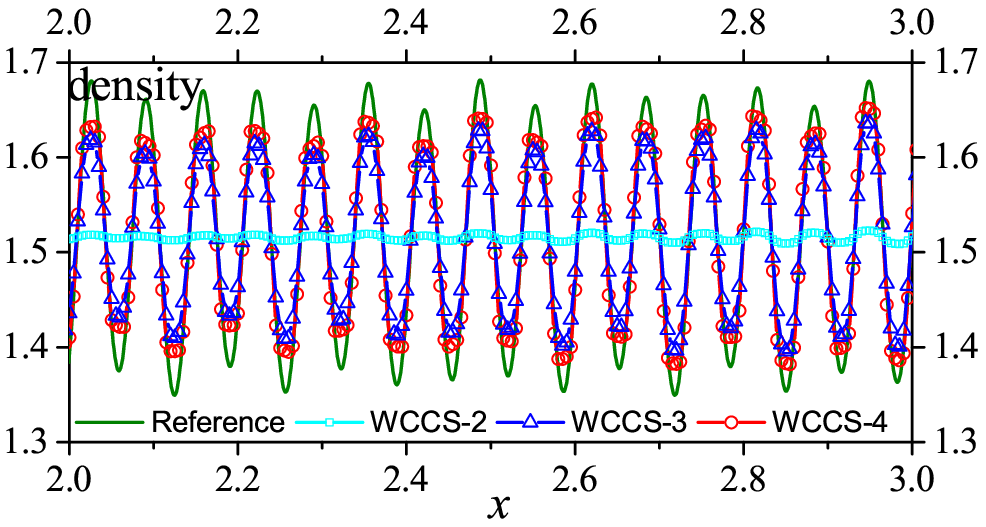}}
  \caption{Density profiles of the modified Shu-Osher problem at $t=5$ calculated by WCCS with $\Delta x=1/200$.}
\label{FIG:Titarev_Toro}
\end{figure}

The second case is the modified Shu-Osher problem \cite{Shu1988EfficientENO} proposed by Titarev and Toro\cite{Titarev2004WENO_FV}.
It depicts the interaction of a high-frequency entropy sine wave with a Mach 1.1 moving shock
which is very suitable for testing the ability of a numerical scheme to capture sophisticated structures.
The computational domain is [-5,5], and the initial condition is given by
\begin{equation*}
  (\rho,u,p)=\begin{cases}
               (1.515695,0.523346,1.805), & \mbox{if } x<-4.5 \\
               (1+0.1\mbox{sin}(20\pi x),0,1), & \mbox{otherwise}.
             \end{cases}
\end{equation*}
Non-reflection boundary conditions are applied on the left and right sides.
Fig. \ref{FIG:Titarev_Toro} shows the density profiles at $t=5$ calculated by WCCS with $\Delta x=1/200$.
The reference solution is calculated by the fifth-order finite-difference WENO-Z scheme \cite{Borges2008WENO_Z} with $\Delta x=1/1000$.
We observe that WCCS-2 damps the wave seriously because it is only first-order at local extrema.
The comparison of the enlarged view demonstrates that WCCS-4 performs the best among the three schemes.

\subsection{Numerical examples for the 2D Euler equations}
The 2D compressible Euler equations can be written as Eq. (\ref{Eq:2DHCLSystem})
with $\mathbf{U}=[\rho,\rho u,\rho v,\rho e]^T$, $\mathbf{F}=[{\rho u,\rho u^2+p,\rho uv,(\rho e+p)u}]^T$,
and $\mathbf{G}=[{\rho v,\rho uv,\rho v^2+p,(\rho e+p)v}]^T$,
where $\rho$, $u$, $v$, $p$, and $e$ denote the density, $x$-velocity, $y$-velocity, pressure, and specific total energy respectively.
We use the prefect gas model and the specific total energy is calculated as $e=\frac{p}{\rho(\gamma-1)}+\frac{1}{2} (u^2+v^2)$.
Here, we set the specific heat ratio $\gamma=1.4$.
When implementing the rotated characteristic decomposition technique,
we use $\phi=\rho e$ in Eq. (\ref{Eq:RotatedNormal}) to determine the rotated direction.
The details of the characteristic information for the Euler equations can be found in Toro's book\cite{Toro2013CFDBook}.
\subsubsection{Isentropic vortex evolution problem}
The computational domain is $[-5, 5]\times[-5, 5]$,  and all boundaries are periodic.
The initial condition is a uniform flow $\rho=1, p=1$ and $(u,v)=(1,1)$
adding a isentropic perturbation expressed as
\begin{equation}\label{Eq:Vortex_Perturbations}
  (\delta u, \delta v)=\frac{\psi}{2\pi}e^{0.5(1-r^2)}(-y, x), \delta T=-\frac{(\gamma-1)\psi^2}{8\gamma\pi^2}e^{(1-r^2)}
   \text{ and }\delta S=0,
\end{equation}
where the entropy $r^2=x^2+y^2$, $S=p/\rho^\gamma$, the temperature $T=p/\rho$, and the vortex strength $\psi=5$.

The exact solution of the problem is the vortex moving at the mean velocity,
so we can use it to test the convergence rate of 2D numerical schemes.
Table \ref{Table:2DVortexConvergence} shows density errors at $t=2$ computed by LCCS and WCCS with different mesh sizes.
All the schemes achieve the expected convergence rates.
The limiter has slight effects on the numerical errors, especially for WCCS-3 and WCCS-4.
The 2nd-order convergence rates of WCCS-2 deserve special attentions,
because we should expect 1st-order convergence for $L_\infty$ of WCCS-2 due to the TVD like limiter.
We found that $L_\infty$ of WCCS-2 indeed converges lower than 2nd-order
without using characteristic decomposition technique.
The reason is not completely clear because of the complexity of the rotated characteristic decomposition.

\begin{table}
  \centering
  \begin{tabular}{|c|c|c|c|c|c|c|c|c|}
    \hline
      & \multicolumn{4}{|c|}{LCCS-2} & \multicolumn{4}{|c|}{WCCS-2} \\
     \hline
      Mesh size       & $L_1$  &Order & $L_\infty$  &Order    & $L_1$  &Order & $L_\infty$  &Order    \\
      \hline
     1/5              &6.51E-4 &-     &1.33E-2      &-        &1.69E-3 &-     &3.77E-2      &- \\

     1/10             &1.15E-4 &2.50  &2.42E-3      &2.46     &2.75E-4 &2.62  &9.63E-3      &1.95 \\

     1/20             &2.17E-5 &2.41  &4.48E-4      &2.43     &4.21E-5 &2.71  &2.03E-3      &2.25 \\

     1/40             &4.79E-6 &2.18  &8.98E-5      &2.32     &6.85E-6 &2.62  &3.74E-4      &2.44 \\
     \hline

     & \multicolumn{4}{|c|}{LCCS-3} & \multicolumn{4}{|c|}{WCCS-3} \\
     \hline
      Mesh size       & $L_1$  &Order & $L_\infty$  &Order    & $L_1$  &Order & $L_\infty$  &Order    \\
     \hline
     1/5              &3.71E-4 &-     &5.73E-3      &-        &4.09E-4 &-     &6.91E-3      &- \\

     1/10             &5.13E-5 &2.85  &7.91E-4      &2.86     &5.26E-5 &2.96  &9.30E-4      &2.89 \\

     1/20             &6.53E-6 &2.97  &1.03E-4      &2.94     &6.90E-6 &2.93  &1.15E-4      &3.02 \\

     1/40             &8.20E-7 &2.99  &1.29E-5      &3.00     &8.20E-7 &3.07  &1.29E-5      &3.16 \\
     \hline

     & \multicolumn{4}{|c|}{LCCS-4} & \multicolumn{4}{|c|}{WCCS-4} \\
     \hline
      Mesh size       & $L_1$  &Order & $L_\infty$  &Order    & $L_1$  &Order & $L_\infty$  &Order    \\
     \hline
     1/5              &5.54E-5 &-     &7.25E-4      &-        &5.70E-5 &-     &7.25E-4      &- \\

     1/10             &2.13E-6 &4.70  &3.41E-5      &4.41     &2.16E-6 &4.72  &3.41E-5      &4.41 \\

     1/20             &7.45E-8 &4.84  &1.20E-6      &4.83     &7.52E-8 &4.84  &1.20E-6      &4.83 \\

     1/40             &4.67E-9 &4.00  &4.26E-8      &4.82     &5.49E-9 &3.78  &5.55E-8      &4.43 \\
     \hline
  \end{tabular}
  \caption{Density errors of the isentropic vortex evolution problem at $t=2$ computed by LCCS and WCCS with different mesh sizes.}
  \label{Table:2DVortexConvergence}
\end{table}

\subsubsection{Two-dimensional Riemann problems}

\begin{figure}
  \centering
  \subfigure[WCCS-2]{
  \label{FIG:RP1_WCCS2}
  \includegraphics[width=5.5 cm]{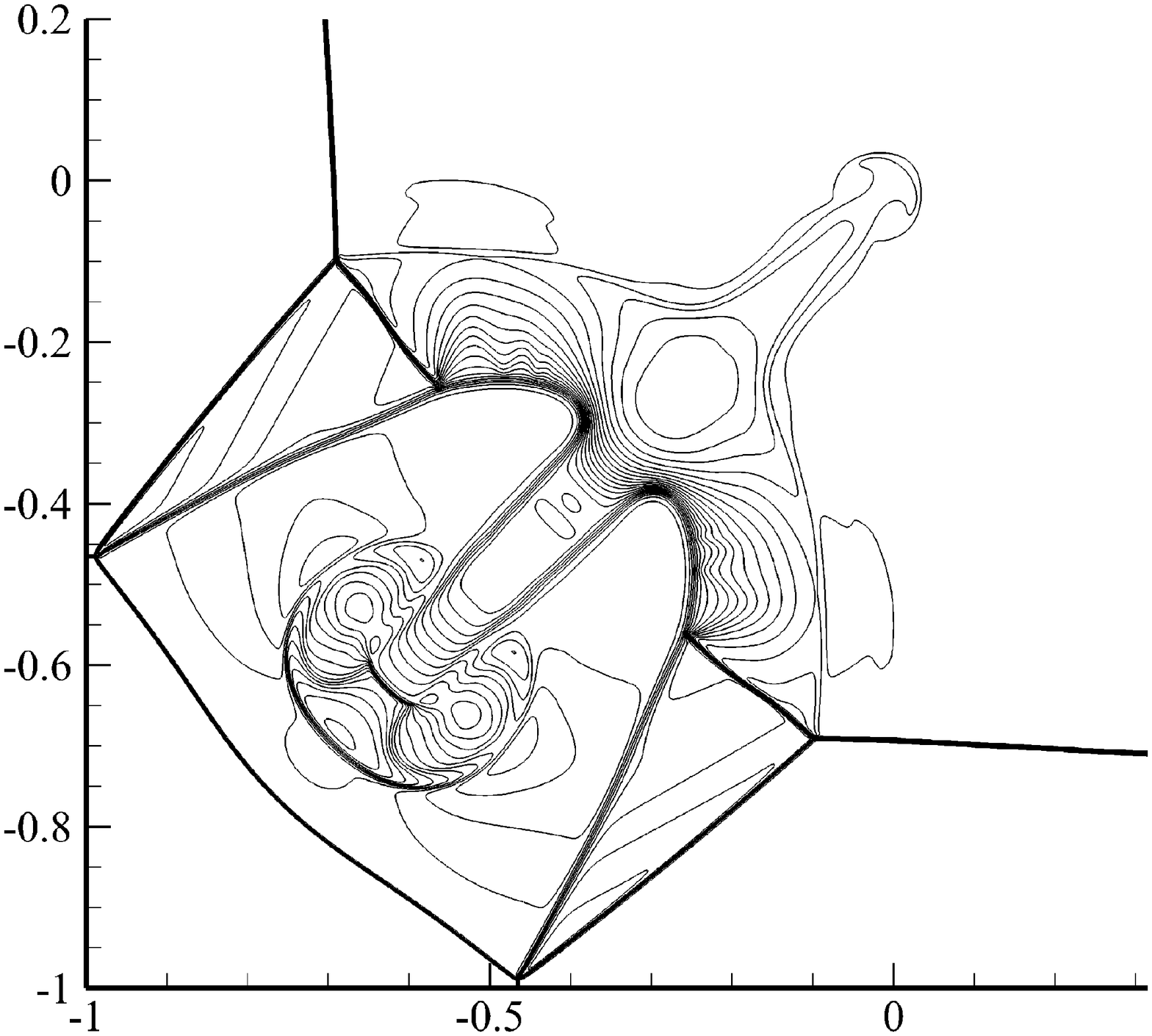}}
  \subfigure[WCCS-3]{
  \label{FIG:RP1_WCCS3}
  \includegraphics[width=5.5 cm]{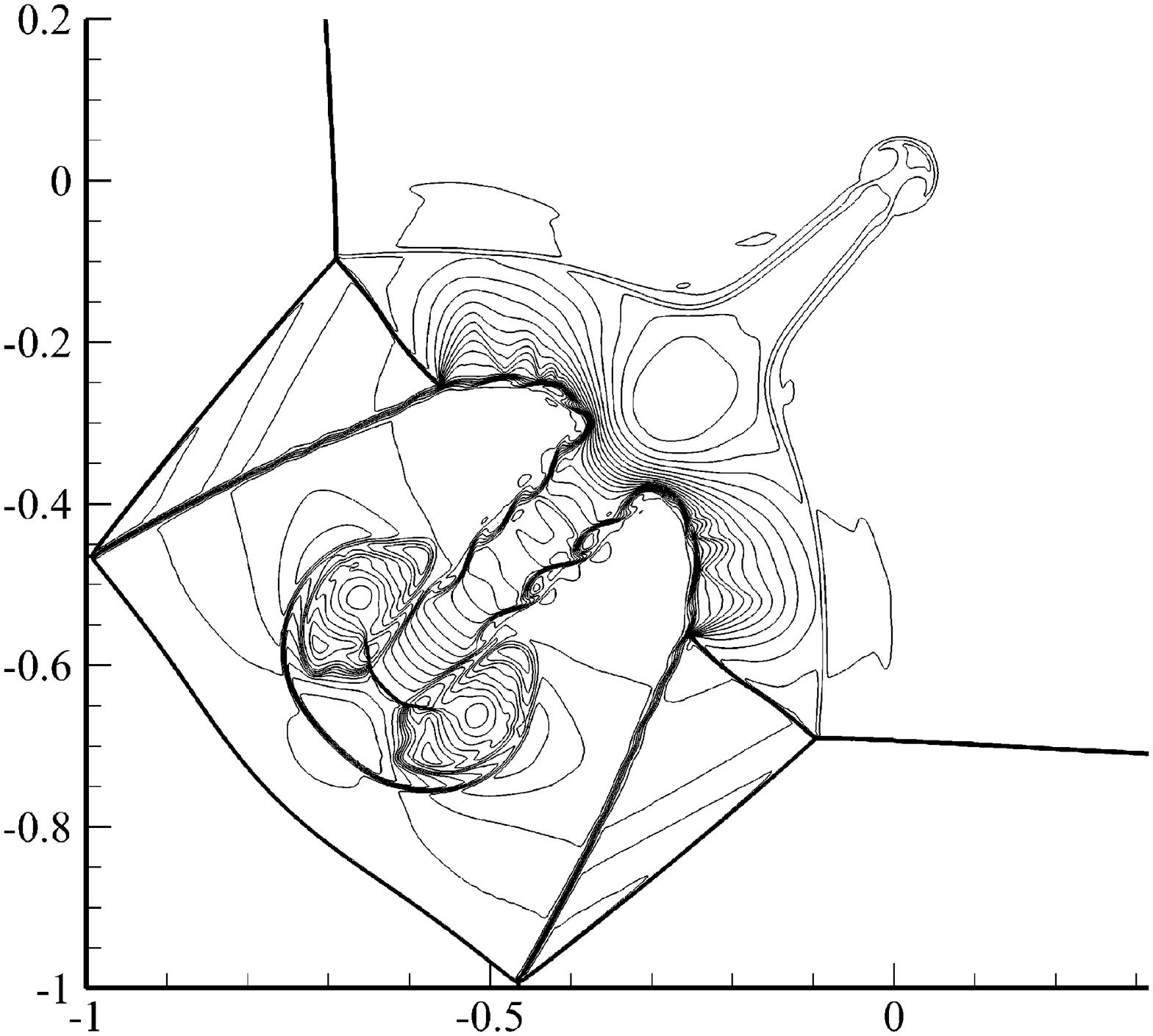}}
  \subfigure[WCCS-4]{
  \label{FIG:RP1_WCCS4}
  \includegraphics[width=5.5 cm]{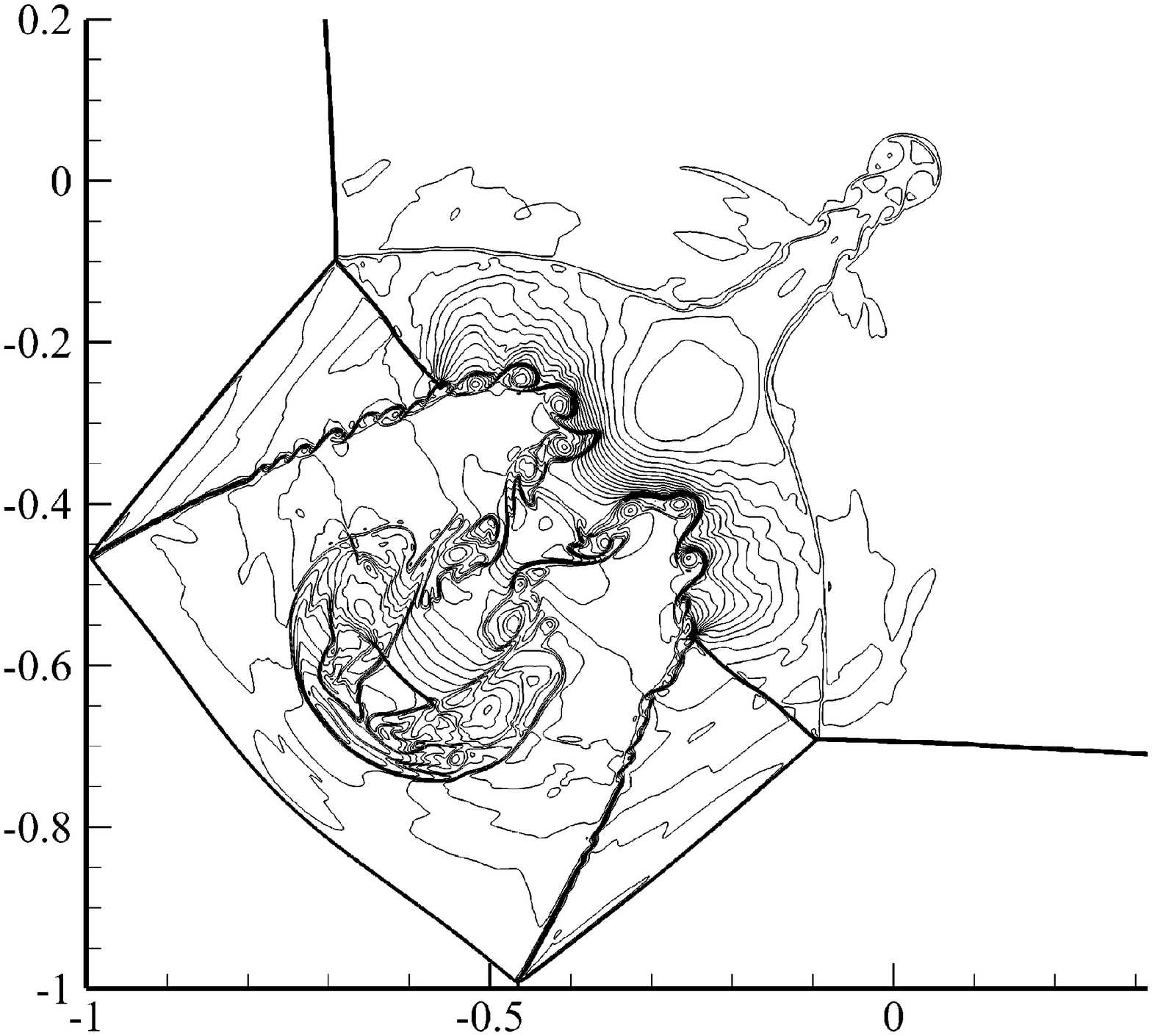}}
  \caption{Density contours of the first Riemann problem at $t=1.1$ calculated by WCCS with $1200\times1200$ cells.}
\label{FIG:RP1_Density_Contour}
\end{figure}

The initial condition of the first Riemann problem is given by
\begin{equation*}
  \left(\rho, u, v, p\right)=
  \begin{cases}
    \left(0.138,1.206,1.206,0.029\right),    & \text{in }[-1, 0]\times[-1, 0],\\
    \left(0.5323,1.206,0,0.3\right),         & \text{in }[-1, 0]\times[0, 1],\\
    \left(1.5,0,0,1.5\right),                & \text{in }[0, 1] \times[0, 1],\\
    \left(0.5323,0,1.206,0.3\right),         & \text{in }[0, 1] \times[-1, 0].
  \end{cases}
\end{equation*}
This configuration depicts the interaction of four normal shocks
which results in two double-Mach reflections and an oblique shock moving along the diagonal of the computational domain.
The density contours at $t=1.1$ calculated by WCCS schemes with $1200\times1200$ cells are shown in Fig. \ref{FIG:RP1_Density_Contour}.
We can see that WCCS-3 begins to trace out the Kelvin–Helmholtz (KH) instabilities which are absent in the results of WCCS-2.
WCCS-4 captures the most details of the vortices among the three schemes.

\begin{figure}
  \centering
  \subfigure[WCCS-2]{
  \label{FIG:RP2_WCCS2}
  \includegraphics[width=5.5 cm]{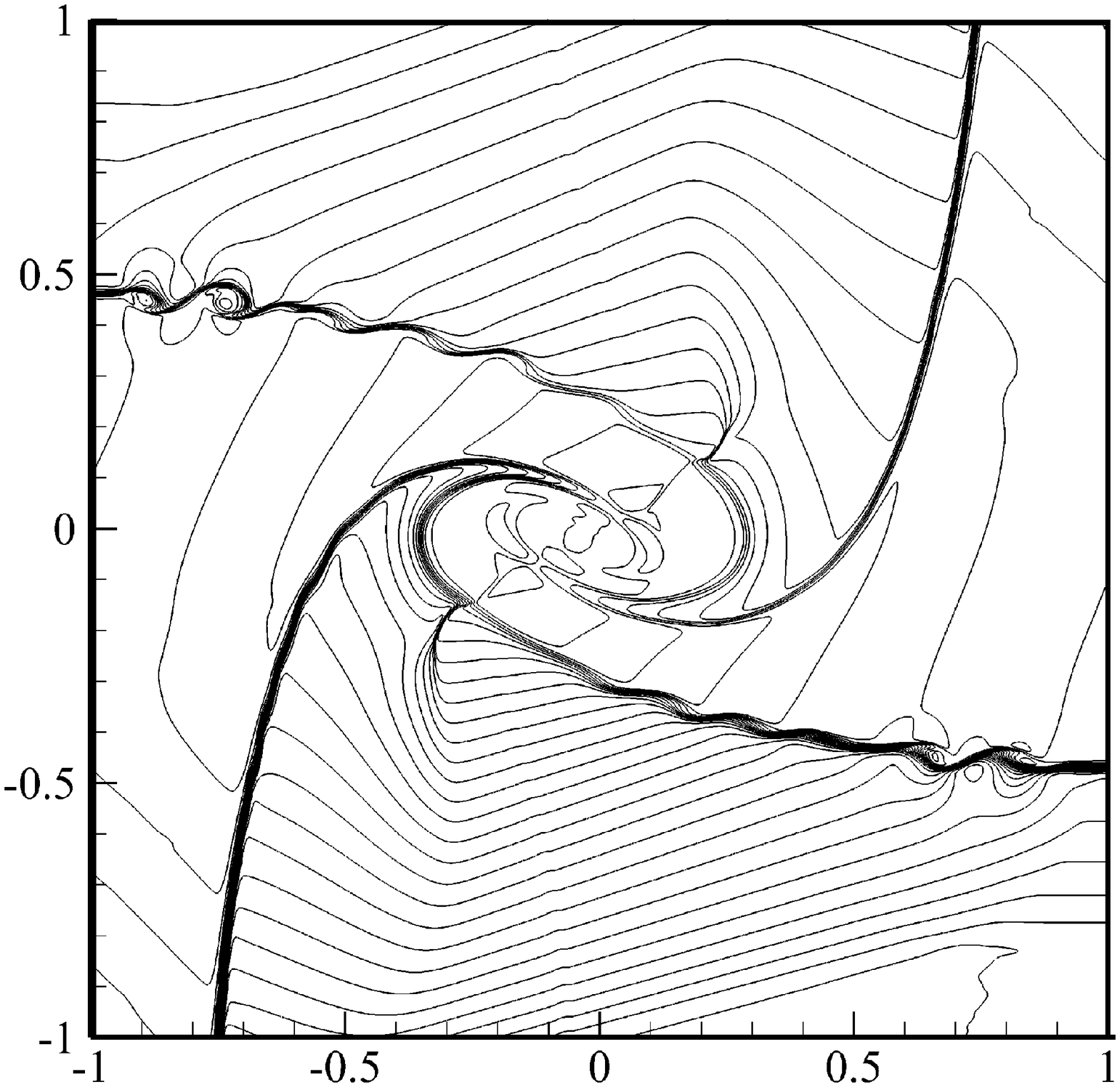}}
  \subfigure[WCCS-3]{
  \label{FIG:RP2_WCCS3}
  \includegraphics[width=5.5 cm]{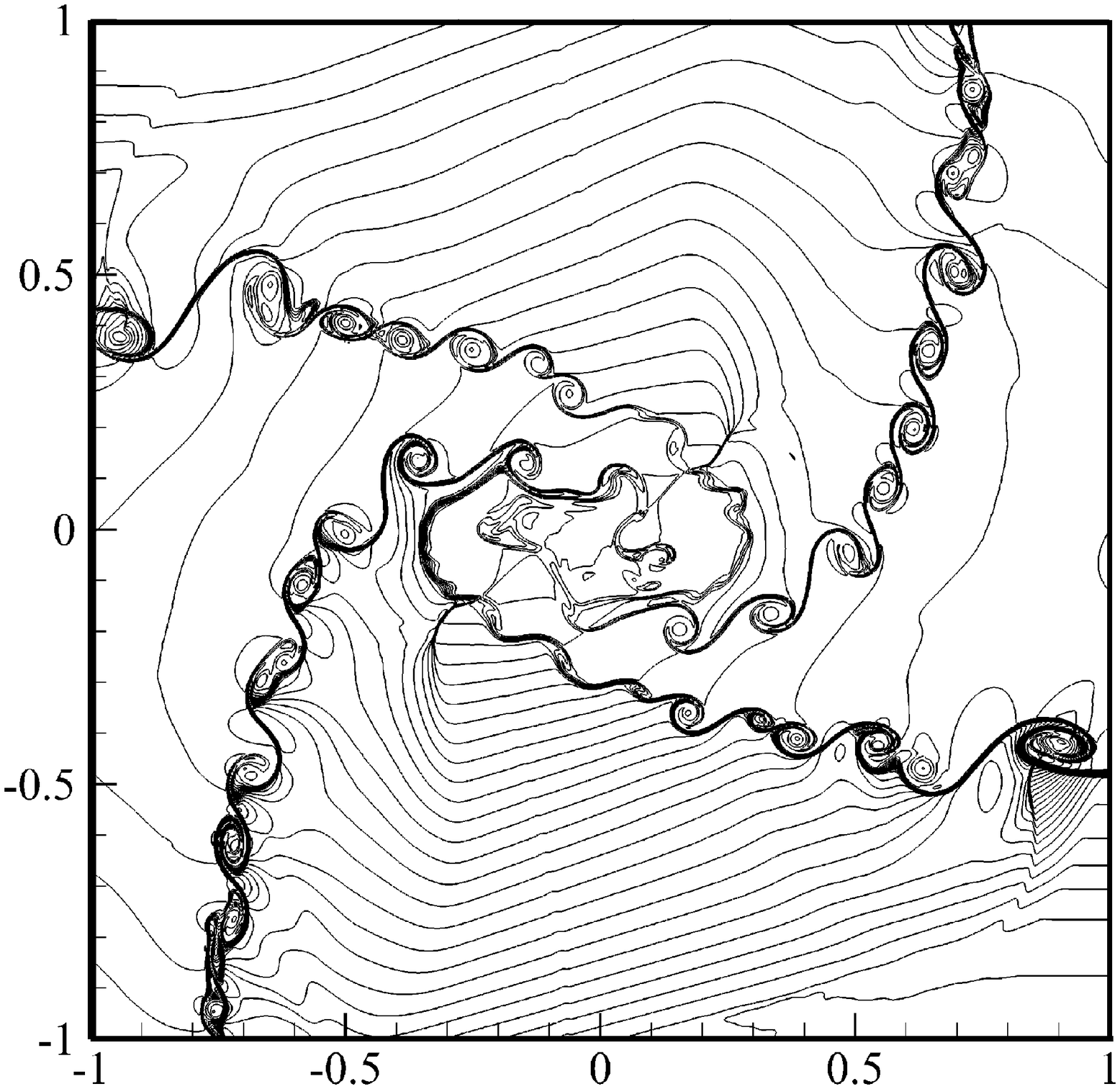}}
  \subfigure[WCCS-4]{
  \label{FIG:RP2_WCCS4}
  \includegraphics[width=5.5 cm]{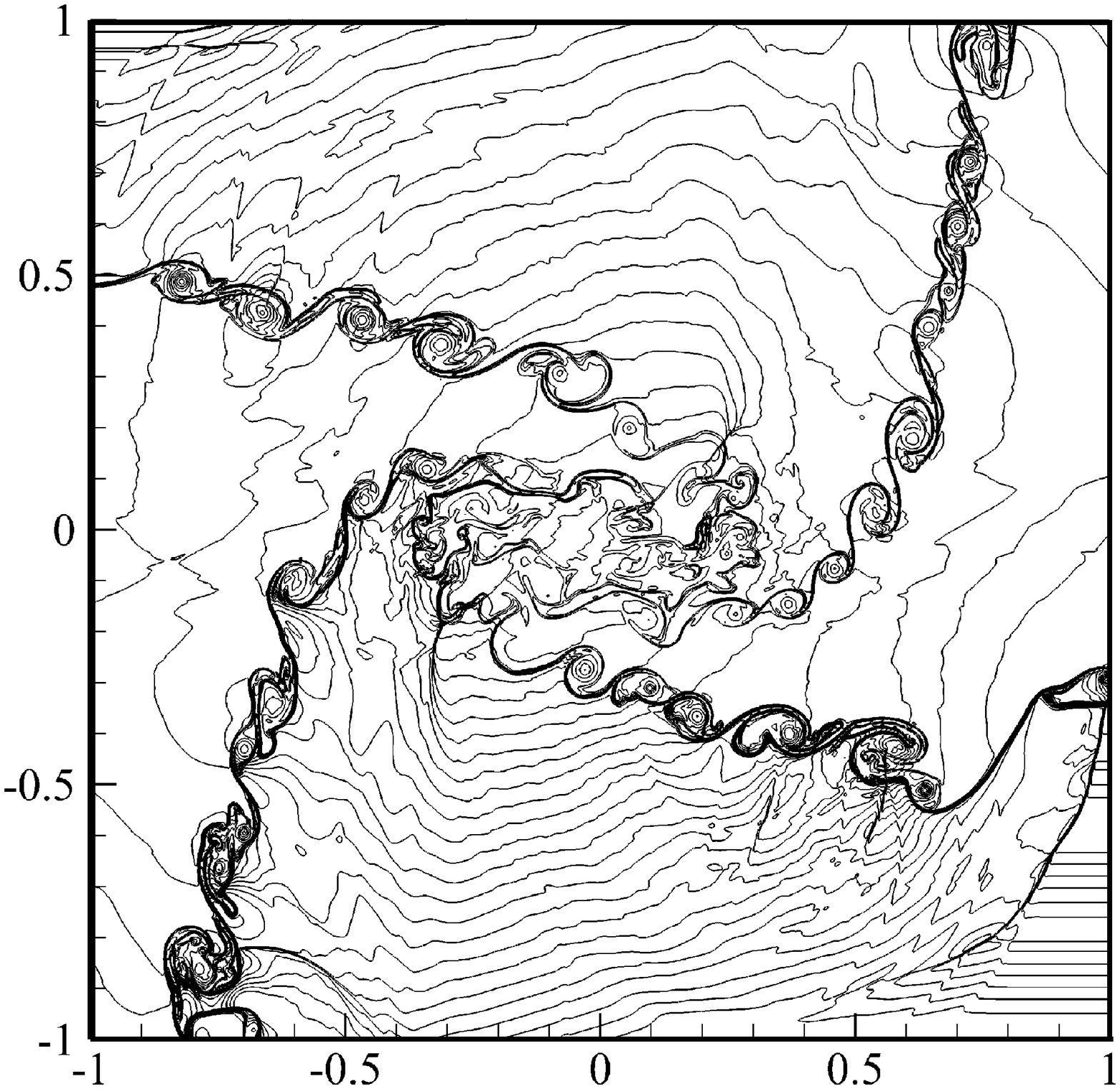}}
  \caption{Density contours of the second Riemann problem at $t=1$ calculated by WCCS with $1200\times1200$ cells.}
\label{FIG:RP2_Density_Contour}
\end{figure}

The second Riemann problem starts off four contact discontinuities which are expressed as
\begin{equation*}
  \left(\rho, u, v, p\right)=
  \begin{cases}
    \left(1,-0.75,0.5,1\right),        & \text{in }[-1, 0]\times[-1, 0],\\
    \left(2,0.75,0.5,1\right),         & \text{in }[-1, 0]\times[0, 1],\\
    \left(1,0.75,-0.5,1\right),        & \text{in }[0, 1] \times[0, 1],\\
    \left(3,-0.75,-0.5,1\right),       & \text{in }[0, 1] \times[-1, 0].
  \end{cases}
\end{equation*}
The density contours at $t=1$ calculated by WCCS schemes with $1200\times1200$ cells are shown in Fig. \ref{FIG:RP2_Density_Contour}.
We observe a large-scale spiral at the center of the domain in all results,
but small-scale vortices at the vicinity of contact discontinuities only appear in the results of WCCS-3 and WCCS-4.

\subsubsection{Double Mach reflection problem}
\begin{figure}
  \centering
  \subfigure[WCCS-2]{
  \label{FIG:DMR_WCCS2}
  \includegraphics[width=10 cm]{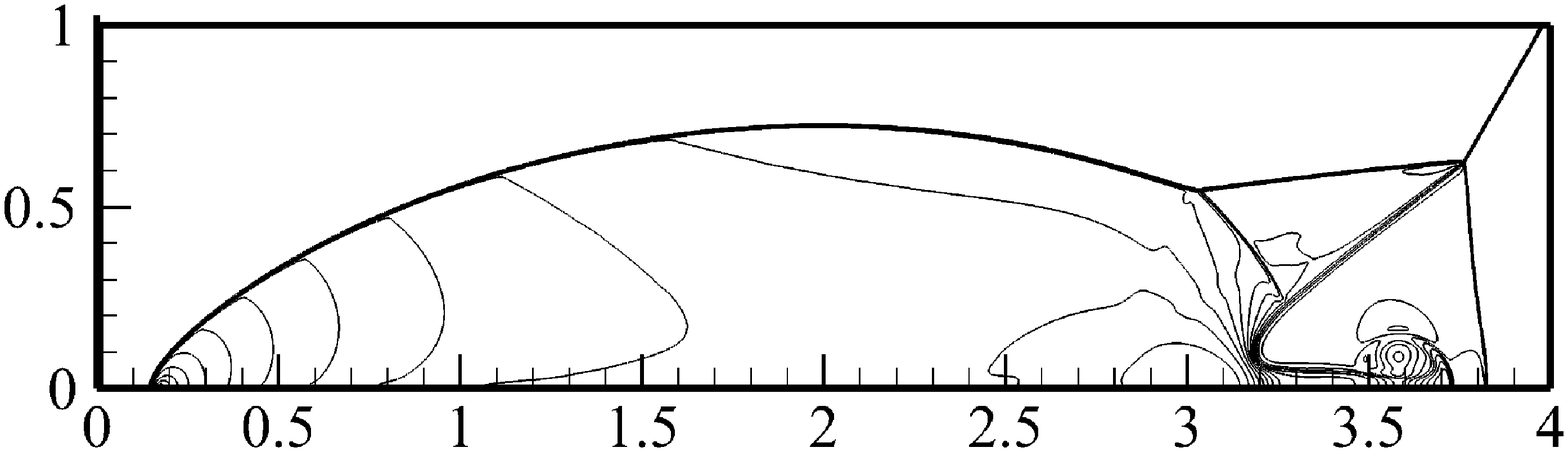}}
  \subfigure[WCCS-3]{
  \label{FIG:DMR_WCCS3}
  \includegraphics[width=10 cm]{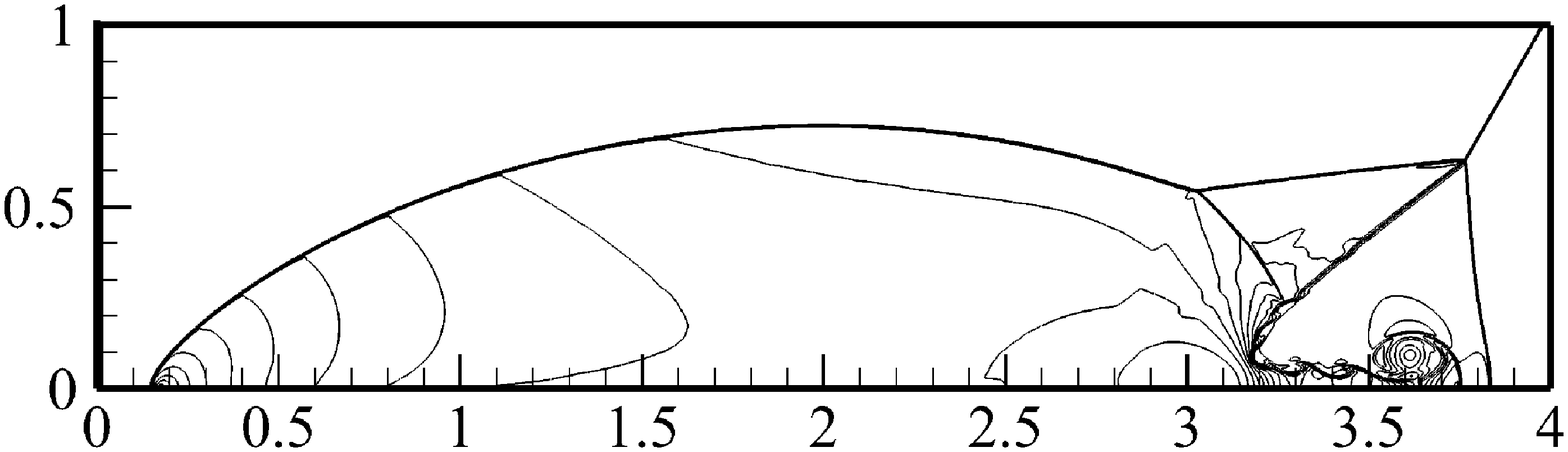}}
  \subfigure[WCCS-4]{
  \label{FIG:DMR_WCCS4}
  \includegraphics[width=10 cm]{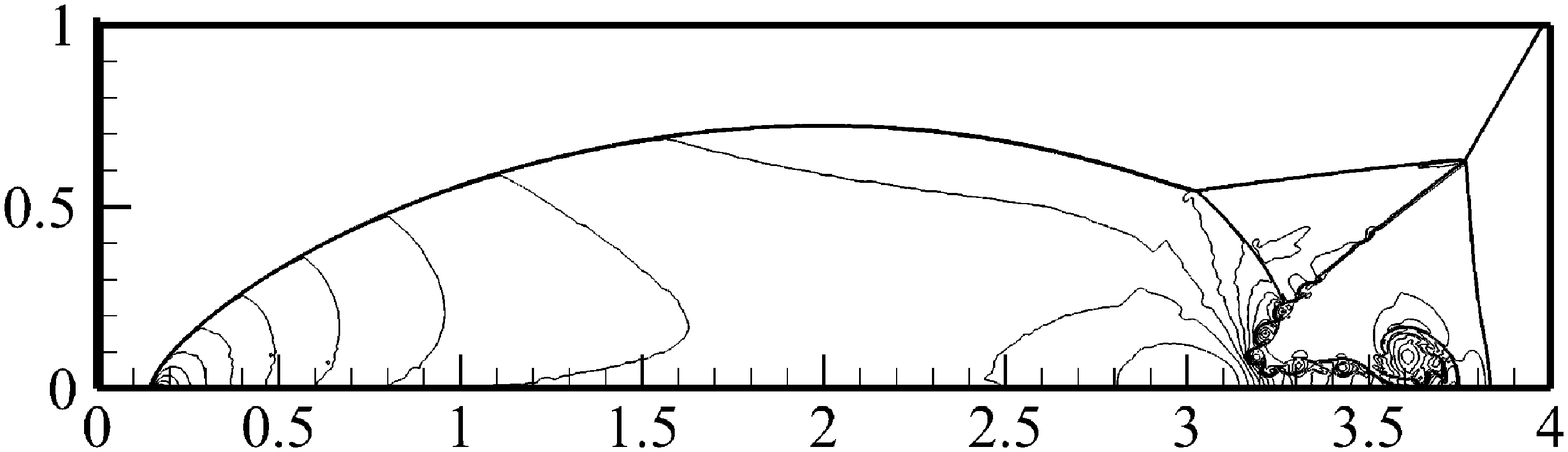}}
  \caption{Density contours of the double Mach reflection problem at $t=0.28$ calculated by WCCS with $1600\times400$ cells.}
\label{FIG:DMR_Density_Contour}
\end{figure}

\begin{figure}
  \centering
  \subfigure[WCCS-2]{
  \label{FIG:DMR_enlarge_WCCS2}
  \includegraphics[width=5.5 cm]{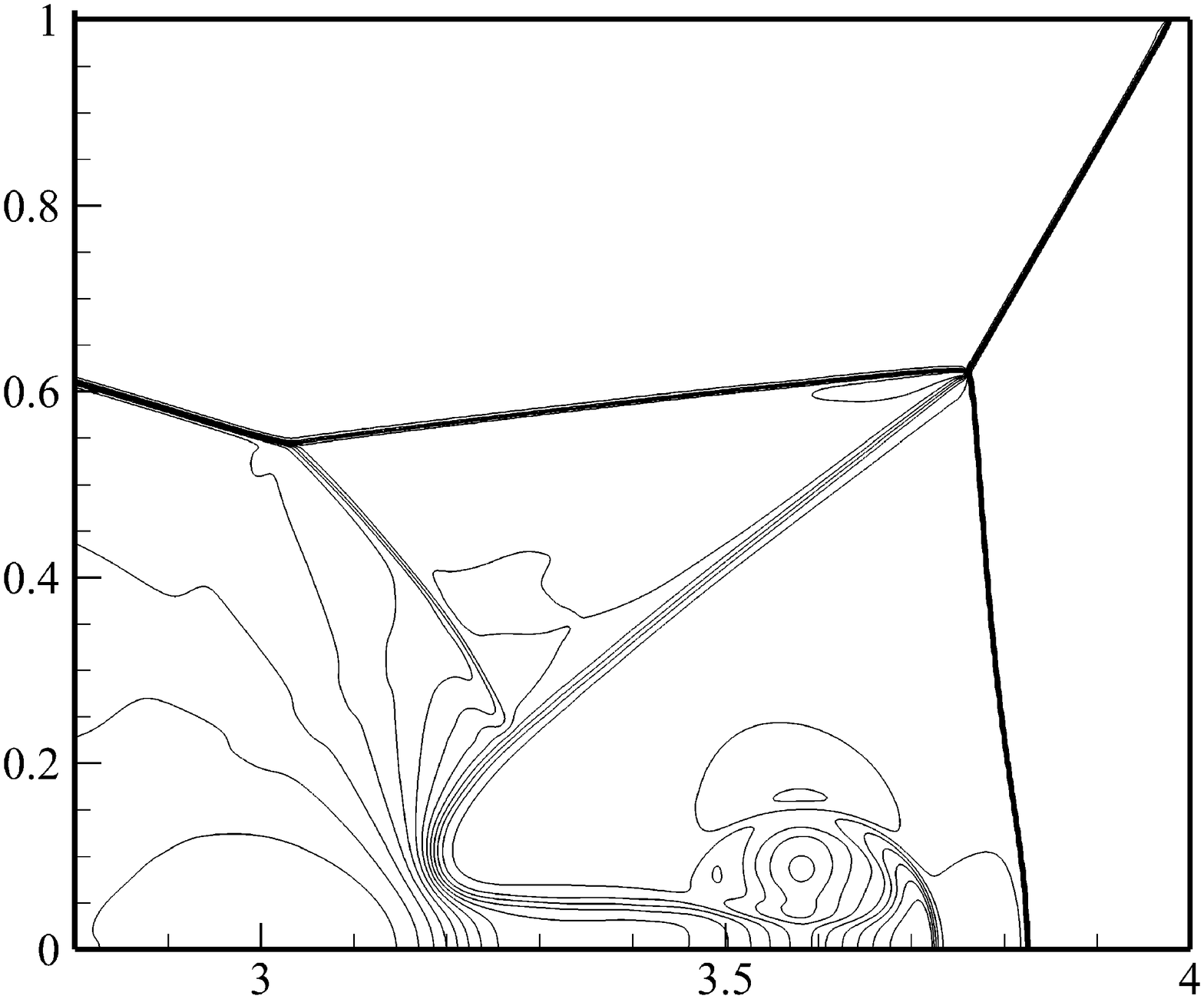}}
  \subfigure[WCCS-3]{
  \label{FIG:DMR_enlarge_WCCS3}
  \includegraphics[width=5.5 cm]{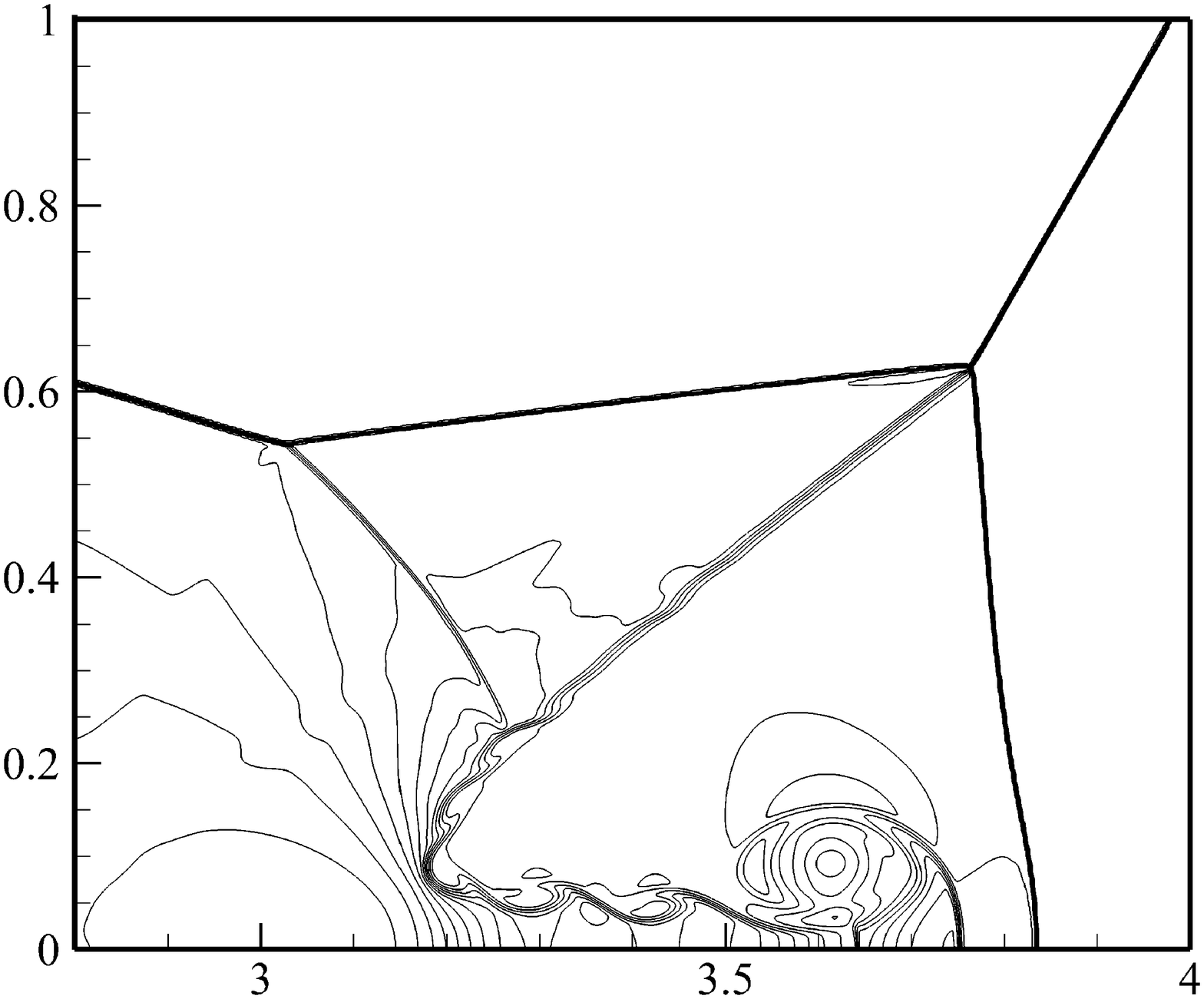}}
  \subfigure[WCCS-4]{
  \label{FIG:DMR_enlarge_WCCS4}
  \includegraphics[width=5.5 cm]{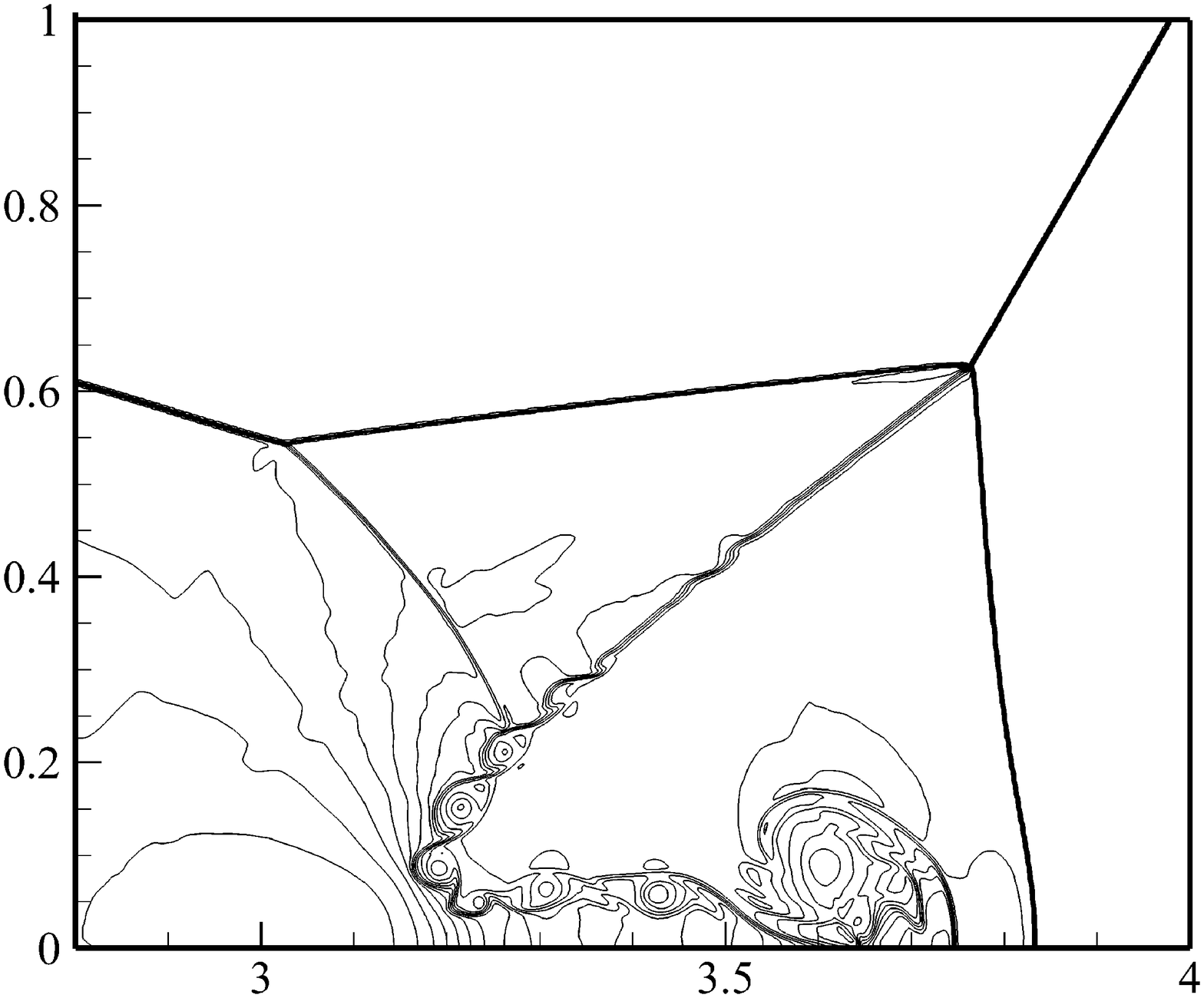}}
  \caption{Enlarged view of the density contours of the double Mach reflection problem at $t=0.28$ calculated by WCCS with $1600\times400$ cells.}
\label{FIG:DMR_Density_Contour_enlarge}
\end{figure}

This problem is a classical test case for high-resolution numerical methods which was originally proposed by Woodward and Colella \cite{Woodward1984JCP}.
At the beginning, the computational domain $[0, 4]\times[0, 1]$ contains an oblique Mach 10 shock
which is inclined at an angle of $60^\circ$ to the horizontal direction.
The pre- and post-shocked states are given by
\begin{equation}\label{Eq:DMR}
\left(\rho, u, v, p\right)=
  \begin{cases}
    \left(1.4, 0, 0, 1\right) & \text{if } x>\frac{1}{6}+\frac{y}{\sqrt{3}}, \\
    \left(8, 8.25sin(60^\circ), -8.25cos(60^\circ), 116.5\right), & \text{otherwise}.
  \end{cases}
\end{equation}
On the left, the post-shocked states are imposed;
On the right, nonreflective boundary conditions are implemented;
On the top, the boundary conditions are determined by the exact motion of the oblique shock;
On the bottom, nonreflective boundary conditions are imposed if $x\le\frac{1}{6}$,
otherwise reflective wall boundary conditions are imposed.
Fig. \ref{FIG:DMR_Density_Contour} shows the entire view of the density contours at $t=0.28$ calculated by WCCS with $1600\times400$ cells.
As we see, all schemes capture the shocks and Mach stems very well.
However, when we zoom in the view as shown in Fig. \ref{FIG:DMR_Density_Contour_enlarge},
we observe that WCCS-3 and WCCS-4 capture the small-scale vortices induced by KH instabilities that are completely dissipated by WCCS-2.

\subsubsection{Single-mode Richtmyer–Meshkov instabilities}

\begin{figure}
  \centering
  \subfigure[WCCS-2]{
  \label{FIG:RMI_WCCS2}
  \includegraphics[height=8 cm]{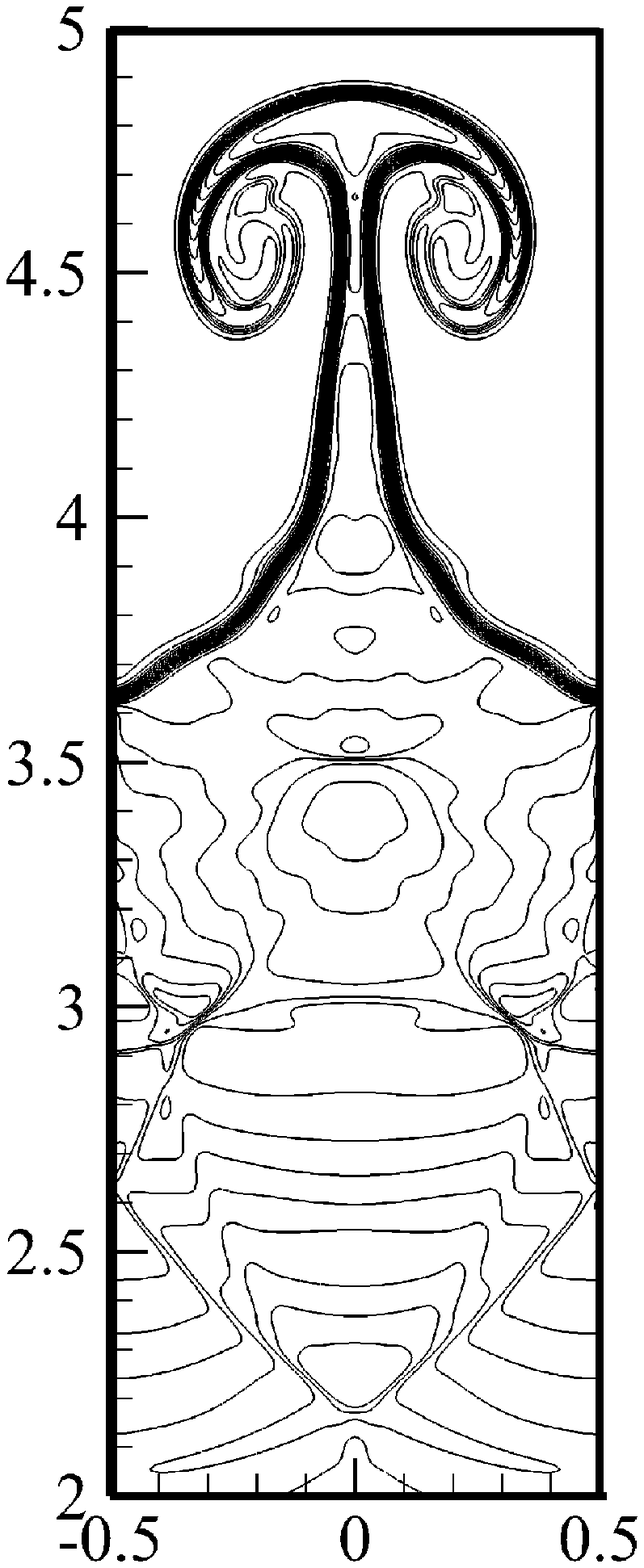}}
  \subfigure[WCCS-3]{
  \label{FIG:RMI_WCCS3}
  \includegraphics[height=8 cm]{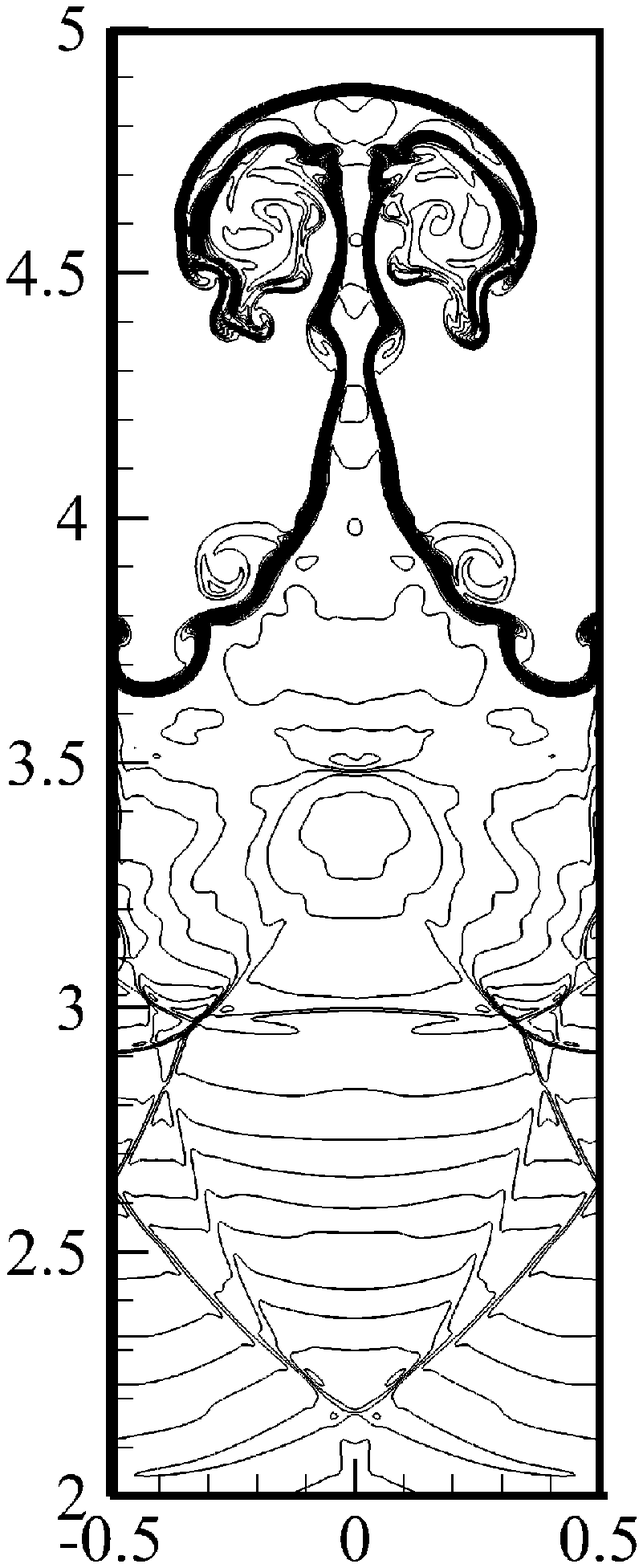}}
  \subfigure[WCCS-4]{
  \label{FIG:RMI_WCCS4}
  \includegraphics[height=8 cm]{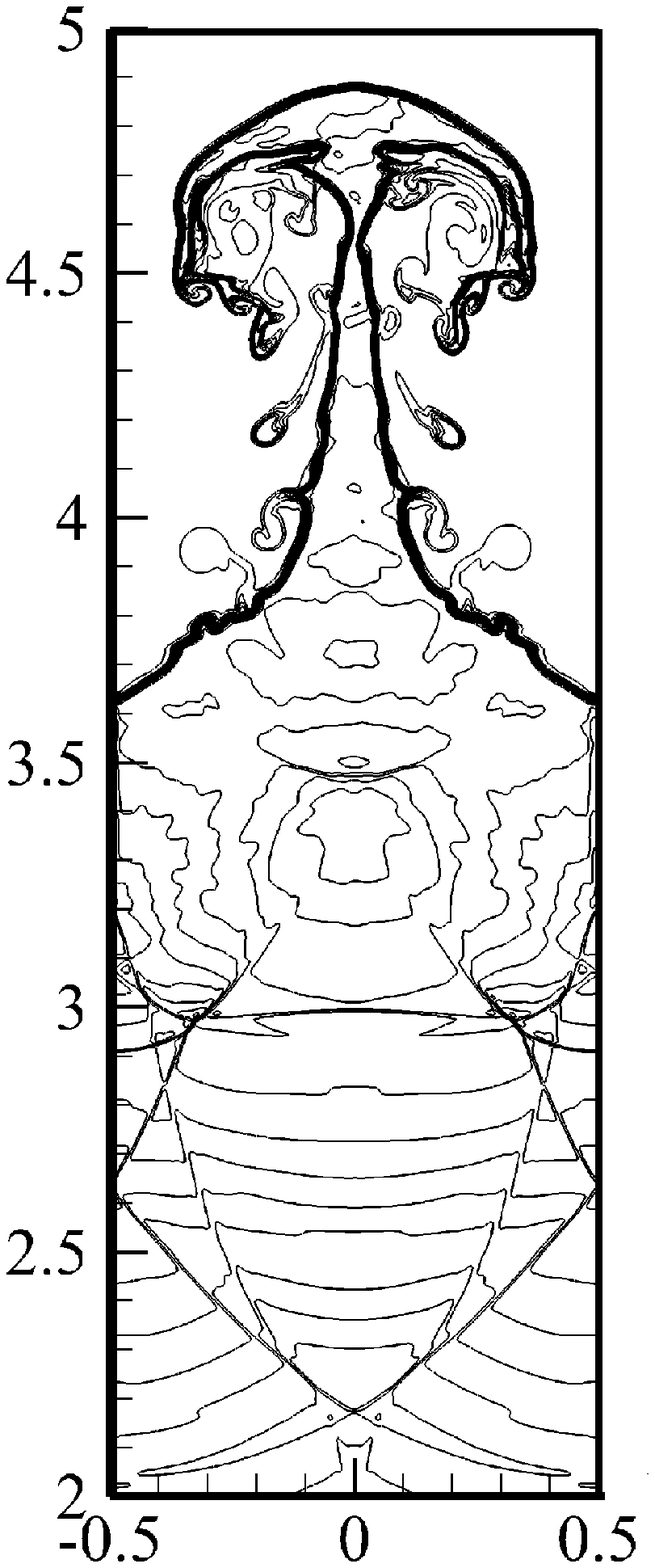}}
  \caption{Density contours of the single-mode Richtmyer–Meshkov instabilities at $t=1.8$ calculated by WCCS with $300\times1500$ cells.}
\label{FIG:RMI_Density_Contour}
\end{figure}

The computational domain is $[-0.5, 0.5]\times[0.0, 5.0]$ with the initial conditions given by
\begin{equation*}
  \left(\rho, u, v, p\right)=
  \begin{cases}
    \left(\rho_1,u_1,v_1,p_1\right),        & \text{if } y\ge 1-0.3\cos(2\pi x),\\
    \left(\rho_2,u_2,v_2,p_2\right),         & \text{else if } y\ge 0.6,\\
    \left(\rho_s,u_s,v_s,p_s\right),         & \text{else.}
  \end{cases}
\end{equation*}
The gas status on two sides of the cosine-shaped contact discontinuity are
$\left(\rho_1,u_1,v_1,p_1\right)=\left(0.1,0,0,1/\gamma\right)$ and
$\left(\rho_2,u_2,v_2,p_2\right)=\left(1.0,0,0,1/\gamma\right)$.
The post-shock status is calculated by the following normal shock relation
\begin{subequations}
  \begin{equation*}
  \frac{\rho_s}{\rho_2}=\frac{0.5\left(\gamma+1\right)M_s^2}{1+0.5\left(\gamma-1\right)M_s^2},
  \end{equation*}
  \begin{equation*}
  u_s=0, \frac{v_s}{a_2}=\frac{2\left(M_s^2-1\right)}{\left(\gamma+1\right)M_s},a_2=\sqrt{\frac{\gamma p}{\rho}},
  \end{equation*}
  \begin{equation*}
  \frac{p_s}{p_2}=\frac{2\gamma M_s^2-\gamma+1}{\gamma+1}.
  \end{equation*}
\end{subequations}
Here, we set the shock Mach number $M_s=2$. The left and right boundaries are periodic,
and the upper and lower boundaries are non-reflective.

After the simulation starts, the shock wave moves towards the cosine-shaped contact discontinuity.
Richtmyer–Meshkov instabilities (RMI) \cite{Brouillette2002RMI_Review} occur
when the contact discontinuity is impulsively accelerated by the impact of the shock wave.
Fig. \ref{FIG:RMI_Density_Contour} shows the density contours at $t=1.8$ calculated by WCCS with $300\times1500$ cells.
Obviously, a higher-order scheme capture a sharper contact discontinuity and more details of the RMI.

\section{Conclusions}
We propose a class of high-order weighted compact central finite-volume-type schemes for solving hyperbolic conservation laws.
The schemes do not need a Riemann solver to calculate the fluxes through the boundaries of a control volume thanks to the space-time staggered mesh.
In addition, the schemes can achieve arbitrarily high-order on a super compact spatial stencil,
and can achieve the same order in space and time with only one explicit time step.
When the order increases, the maximum CFL number for the numerical stability becomes smaller,
but it is still at a relatively high level.
A series of numerical examples demonstrate that the proposed schemes can capture discontinuities without spurious oscillations
while maintaining high order for smooth solutions.
High-order schemes perform much better than the second-order scheme for the problems evolving small-scale structures.

\section*{Acknowledgement}
H. S. would like to acknowledge the financial support of the National Natural Science Foundation of China (Contract No. 11901602).

\begin{appendix}
\section{}\label{AppendixA}
The coefficient matrices $\mathbf{M}_1$ and $\mathbf{M}_2$ in Eq. (\ref{Eq:1DSchemeMatrixForm})
for the three linear solutions (\ref{Eq:1D_U012}) of 1D 2nd-, 3rd-, and 4th-order WCC schemes are derived as follows:
\begin{itemize}
  \item For the 1D 2nd-order WCC scheme, the corresponding coefficient matrices of the three linear solutions are
  \begin{subequations}\label{Eq:2ndOrderCoeMatrix}
    \begin{equation}
    \mathbf{M}_{01}=\begin{bmatrix}
                      \frac{1}{2}+\nu & \frac{1}{8}-\frac{\nu^2}{2} \\
                      -1 & \nu
                    \end{bmatrix},
    \mathbf{M}_{02}=\begin{bmatrix}
                      \frac{1}{2}-\nu & -\frac{1}{8}+\frac{\nu^2}{2} \\
                      1 & -\nu
                    \end{bmatrix},
    \end{equation}

    \begin{equation}
    \mathbf{M}_{11}=\begin{bmatrix}
                      \frac{1}{2}+\nu & \frac{1}{8}-\frac{\nu^2}{2} \\
                      -1+2\nu & \frac{1}{4}+2\nu-\nu^2
                    \end{bmatrix},
    \mathbf{M}_{12}=\begin{bmatrix}
                      \frac{1}{2}-\nu & -\frac{1}{8}+\frac{\nu^2}{2} \\
                      1-2\nu & -\frac{1}{4}+\nu^2
                    \end{bmatrix},
    \end{equation}

    \begin{equation}
    \mathbf{M}_{21}=\begin{bmatrix}
                      \frac{1}{2}+\nu & \frac{1}{8}-\frac{\nu^2}{2} \\
                      -1-2\nu & \nu^2-\frac{1}{4}
                    \end{bmatrix},
    \mathbf{M}_{22}=\begin{bmatrix}
                      \frac{1}{2}-\nu & -\frac{1}{8}+\frac{\nu^2}{2} \\
                      1+2\nu & \frac{1}{4}-2\nu-\nu^2
                    \end{bmatrix},
    \end{equation}
  \end{subequations}
  where $\nu=a\frac{\Delta t}{\Delta x}$, and the first subscript of $\mathbf{M}$
  corresponds to the subscript of $u$ in Eqs. (\ref{Eq:1D_U012}) and (\ref{Eq:1DLimitedPolynomial}).
  When $\nu\in[0,0.5]$, the stability conditions Eq. (\ref{Eq:LinearStabilityCondition}) of the three 2nd-order linear solutions are satisfied,
  so the linear stability condition of the 2nd-order WCC is $\nu\in[0,0.5]$.

  \item For the 1D 3rd-order WCC scheme, the corresponding coefficient matrices of the three linear solutions are
  \begin{subequations}\label{Eq:2ndOrderCoeMatrix}
    \begin{equation}
    \mathbf{M}_{01}=\begin{bmatrix}
                      \frac{1}{2}+\nu & \frac{1}{6}-\frac{\nu^2}{2} & \frac{1}{48}-\frac{\nu}{24}+\frac{\nu^3}{6}\\
                      -1              & \nu                         &-\frac{\nu^2}{2}\\
                      0               & -1                          & \nu
                    \end{bmatrix},
    \end{equation}
    \begin{equation}
    \mathbf{M}_{02}=\begin{bmatrix}
                      \frac{1}{2}-\nu & -\frac{1}{6}+\frac{\nu^2}{2} & \frac{1}{48}+\frac{\nu}{24}-\frac{\nu^3}{6}\\
                      1              & -\nu                         &\frac{\nu^2}{2}\\
                      0               & 1                          & -\nu
                    \end{bmatrix},
    \end{equation}

    \begin{equation}
    \mathbf{M}_{11}=\begin{bmatrix}
                      \frac{1}{2}+\nu & \frac{1}{8}-\frac{\nu^2}{2} & \frac{1}{48}+\frac{\nu^3}{6}\\
                      -1+2\nu         & \frac{1}{4}+2\nu-\nu^2      & \frac{1}{24}-\nu^2+\frac{\nu^3}{3}\\
                      0               & 0                           & 0
                    \end{bmatrix},
    \end{equation}
    \begin{equation}
    \mathbf{M}_{12}=\begin{bmatrix}
                      \frac{1}{2}-\nu & -\frac{1}{8}+\frac{\nu^2}{2} & \frac{1}{48}-\frac{\nu^3}{6}\\
                      1-2\nu & -\frac{1}{4}+\nu^2 & \frac{1}{24}-\frac{\nu^3}{3}\\
                      0               & 0                          & 0
                    \end{bmatrix},
    \end{equation}

    \begin{equation}
    \mathbf{M}_{21}=\begin{bmatrix}
                      \frac{1}{2}+\nu & \frac{1}{8}-\frac{\nu^2}{2} & \frac{1}{48}+\frac{\nu^3}{6}\\
                      -1-2\nu         & -\frac{1}{4}+\nu^2      & -\frac{1}{24}-\frac{\nu^3}{3}\\
                      0               & 0                           & 0
                    \end{bmatrix},
    \end{equation}
    \begin{equation}
    \mathbf{M}_{22}=\begin{bmatrix}
                      \frac{1}{2}-\nu & -\frac{1}{8}+\frac{\nu^2}{2} & \frac{1}{48}-\frac{\nu^3}{6}\\
                      1+2\nu & \frac{1}{4}-2\nu-\nu^2 & -\frac{1}{24}+\nu^2+\frac{\nu^3}{3}\\
                      0               & 0                          & 0
                    \end{bmatrix}.
    \end{equation}
  \end{subequations}
  When $\nu\in[0,0.384]$, $u_0$ is stable. When $\nu\in[0,0.5]$, $u_1$ and $u_2$ are stable.

  \item For the 1D 4th-order WCC scheme, the corresponding coefficient matrices of the three linear solutions are
  \begin{subequations}\label{Eq:2ndOrderCoeMatrix}
    \begin{equation}
    \mathbf{M}_{01}=\begin{bmatrix}
                      \frac{1}{2}+\nu & \frac{1}{6}-\frac{\nu^2}{2} & \frac{1}{48}-\frac{\nu}{24}+\frac{\nu^3}{6} & \frac{1}{384}+\frac{\nu^2}{48}-\frac{\nu^4}{24}\\
                      -1              & \nu                         & -\frac{\nu^2}{2}+\frac{1}{24} & -\frac{\nu}{24}+\frac{\nu^3}{6}\\
                      0               & -1                          & \nu & -\frac{\nu^2}{2}\\
                      0               & 0                           & -1  & \nu
                    \end{bmatrix},
    \end{equation}
    \begin{equation}
    \mathbf{M}_{02}=\begin{bmatrix}
                      \frac{1}{2}-\nu & -\frac{1}{6}+\frac{\nu^2}{2} & \frac{1}{48}+\frac{\nu}{24}-\frac{\nu^3}{6} & -\frac{1}{384}-\frac{\nu^2}{48}+\frac{\nu^4}{24}\\
                      1              & -\nu                         & \frac{\nu^2}{2}-\frac{1}{24} & \frac{\nu}{24}-\frac{\nu^3}{6}\\
                      0               & 1                          & -\nu    & \frac{\nu^2}{2}\\
                      0               & 0                           & 1  & -\nu
                    \end{bmatrix},
    \end{equation}

    \begin{equation}
    \mathbf{M}_{11}=\begin{bmatrix}
                      \frac{1}{2}+\nu & \frac{1}{8}-\frac{\nu^2}{2} & \frac{1}{48}+\frac{\nu^3}{6}       & \frac{1}{384}-\frac{\nu^4}{24}\\
                      -1+2v           & \frac{1}{4}+2\nu-\nu^2      & \frac{1}{24}-\nu^2+\frac{\nu^3}{3} & \frac{1}{192}+\frac{\nu^3}{3}-\frac{\nu^4}{12}\\
                      0               & 0                          & 0 & 0\\
                      0               & 0                           & 0  & 0
                    \end{bmatrix},
    \end{equation}
    \begin{equation}
    \mathbf{M}_{12}=\begin{bmatrix}
                      \frac{1}{2}-\nu & -\frac{1}{8}+\frac{\nu^2}{2} & \frac{1}{48}-\frac{\nu^3}{6} & -\frac{1}{384}+\frac{\nu^4}{24}\\
                      1-2\nu          & -\frac{1}{4}+\nu^2          & \frac{1}{24}-\frac{\nu^3}{3} & -\frac{1}{192}+\frac{\nu^4}{12}\\
                      0               & 0                          & 0    & 0\\
                      0               & 0                           & 0  & 0
                    \end{bmatrix},
    \end{equation}

    \begin{equation}
    \mathbf{M}_{21}=\begin{bmatrix}
                      \frac{1}{2}+\nu & \frac{1}{8}-\frac{\nu^2}{2} & \frac{1}{48}+\frac{\nu^3}{6} & \frac{1}{384}-\frac{\nu^4}{24}\\
                      -1-2\nu         & -\frac{1}{4}+\nu^2          & -\frac{1}{24}-\frac{\nu^3}{3} & -\frac{1}{192}+\frac{\nu^4}{12}\\
                      0               & 0                          & 0 & 0\\
                      0               & 0                           & 0  & 0
                    \end{bmatrix},
    \end{equation}
    \begin{equation}
    \mathbf{M}_{22}=\begin{bmatrix}
                      \frac{1}{2}-\nu & -\frac{1}{8}+\frac{\nu^2}{2} & \frac{1}{48}-\frac{\nu^3}{6} & -\frac{1}{384}+\frac{\nu^4}{24}\\
                      1+2\nu          & \frac{1}{4}-2\nu-\nu^2       & -\frac{1}{24}+\nu^2+\frac{\nu^3}{3} & \frac{1}{192}-\frac{\nu^3}{3}-\frac{\nu^4}{12}\\
                      0               & 0                          & 0    & 0\\
                      0               & 0                           & 0  & 0
                    \end{bmatrix}.
    \end{equation}
  \end{subequations}
  When $\nu\in[0,0.304]$, $u_0$ is stable. When $\nu\in[0,0.5]$, $u_1$ and $u_2$ are stable.
\end{itemize}
\section{}\label{AppendixB}
All the terms in Eq. (\ref{Eq:2DIntegral}) are calculated as
\begin{subequations}\label{Eq:2DIntegralTerms}
  \begin{equation}\label{SubEq:2DCellAVG}
  \begin{split}
     \bar{u}_{i+1/2,j+1/2}^{n+1/2}&=\frac{1}{\Delta x\Delta y}\int_{x_{i}}^{x_{i+1}}\int_{y_{j}}^{y_{j+1}}u(x,y,t^{n+1/2})dxdy \\
       &=\sum_{l=0}^{P}\sum_{k=0}^{P-l} \frac{1+(-1)^k}{(k+1)!2^{(k+1)}}\frac{1+(-1)^l}{(l+1)!2^{(l+1)}}(u_{kx,ly})_{i+1/2,j+1/2}^{n+1/2},
  \end{split}
  \end{equation}

  \begin{equation}\label{SubEq:2DCellAVGRU}
  \begin{split}
     (\bar{u}_{RU})_{i,j}^{n}&=\frac{4}{\Delta x\Delta y}\int_{x_{i}}^{x_{i+1/2}}\int_{y_{j}}^{y_{j+1/2}}u(x,y,t^{n}) \\
       &=\sum_{l=0}^{P}\sum_{k=0}^{P-l} \frac{1}{(k+1)!2^k}\frac{1}{(l+1)!2^l}(u_{kx,ly})_{i,j}^{n},
  \end{split}
  \end{equation}

  \begin{equation}\label{SubEq:2DCellAVGRD}
  \begin{split}
     (\bar{u}_{RD})_{i,j+1}^{n}&=\frac{4}{\Delta x\Delta y}\int_{x_{i}}^{x_{i+1/2}}\int_{y_{j+1/2}}^{y_{j+1}}u(x,y,t^{n}) \\
       &=\sum_{l=0}^{P}\sum_{k=0}^{P-l} \frac{1}{(k+1)!2^k}\frac{(-1)^l}{(l+1)!2^l}(u_{kx,ly})_{i,j+1}^{n},
  \end{split}
  \end{equation}

  \begin{equation}\label{SubEq:2DCellAVGLU}
  \begin{split}
     (\bar{u}_{LU})_{i+1,j}^{n}&=\frac{4}{\Delta x\Delta y}\int_{x_{i+1/2}}^{x_{i+1}}\int_{y_{j}}^{y_{j+1/2}}u(x,y,t^{n}) \\
       &=\sum_{l=0}^{P}\sum_{k=0}^{P-l} \frac{(-1)^k}{(k+1)!2^k}\frac{1}{(l+1)!2^l}(u_{kx,ly})_{i+1,j}^{n},
  \end{split}
  \end{equation}

  \begin{equation}\label{SubEq:2DCellAVGLD}
  \begin{split}
     (\bar{u}_{LD})_{i+1,j+1}^{n}&=\frac{4}{\Delta x\Delta y}\int_{x_{i+1/2}}^{x_{i+1}}\int_{y_{j+1/2}}^{y_{j+1}}u(x,y,t^{n}) \\
       &=\sum_{l=0}^{P}\sum_{k=0}^{P-l} \frac{(-1)^k}{(k+1)!2^k}\frac{(-1)^l}{(l+1)!2^l}(u_{kx,ly})_{i+1,j+1}^{n},
  \end{split}
  \end{equation}

  \begin{equation}\label{SubEq:2DFluxFD}
  \begin{split}
     (\bar{f}_D)_{i,j}^n&=\frac{2}{\Delta y\Delta t}\int_{y_{j-1/2}}^{y_j}\int_{t^n}^{t^{n+1/2}}f(x_{i},y,t)dydt \\
       &=\sum_{l=0}^{P}\sum_{k=0}^{P-l} \frac{(-1)^k}{(k+1)!2^k}\frac{1}{(l+1)!}(f_{ky,lt})_{i,j}^{n},
  \end{split}
  \end{equation}

  \begin{equation}\label{SubEq:2DFluxFU}
  \begin{split}
     (\bar{f}_U)_{i,j}^n&=\frac{2}{\Delta y\Delta t}\int_{y_j}^{y_{j+1/2}}\int_{t^n}^{t^{n+1/2}}f(x_{i},y,t)dydt \\
       &=\sum_{l=0}^{P}\sum_{k=0}^{P-l} \frac{1}{(k+1)!2^k}\frac{1}{(l+1)!}(f_{ky,lt})_{i,j}^{n},
  \end{split}
  \end{equation}

  \begin{equation}\label{SubEq:2DFluxGL}
  \begin{split}
     (\bar{g}_L)_{i,j}^n&=\frac{2}{\Delta x\Delta t}\int_{x_{i-1/2}}^{x_i}\int_{t^n}^{t^{n+1/2}}g(x,y_j,t)dxdt \\
       &=\sum_{l=0}^{P}\sum_{k=0}^{P-l} \frac{(-1)^k}{(k+1)!2^k}\frac{1}{(l+1)!}(g_{kx,lt})_{i,j}^{n},
  \end{split}
  \end{equation}

  \begin{equation}\label{SubEq:2DFluxGR}
  \begin{split}
     (\bar{g}_R)_{i,j}^n&=\frac{2}{\Delta x\Delta t}\int_{x_i}^{x_{i+1/2}}\int_{t^n}^{t^{n+1/2}}g(x,y_j,t)dxdt \\
       &=\sum_{l=0}^{P}\sum_{k=0}^{P-l} \frac{1}{(k+1)!2^k}\frac{1}{(l+1)!}(g_{kx,lt})_{i,j}^{n}.
  \end{split}
  \end{equation}
\end{subequations}
\end{appendix}
\bibliography{mybibfile}

\end{document}